\documentclass{amsart}
\usepackage{amsthm,amssymb,amscd}
\usepackage[pagebackref]{hyperref}

\newcommand{\N}{{\mathbb N}}
\newcommand{\Q}{\mathbb Q}
\newcommand{\F}{{\mathbb F}}
\newcommand{\Z}{\mathbb Z}
\newcommand{\OO}{\mathcal O}
\newcommand{\Cl}{\operatorname{Cl}}

\newcommand{\Gal}{\operatorname{Gal}}
\newcommand{\ft}{\mathfrak 2}
\newcommand{\fr}{\mathfrak 3}
\newcommand{\fv}{\mathfrak 5}
\newcommand{\fa}{\mathfrak a}
\newcommand{\fA}{\mathfrak A}
\newcommand{\fb}{\mathfrak b}
\newcommand{\fB}{\mathfrak B}
\newcommand{\fc}{\mathfrak c}
\newcommand{\fl}{\mathfrak l}
\newcommand{\fm}{\mathfrak m}
\newcommand{\fM}{\mathfrak M}
\newcommand{\fp}{{\mathfrak p}}
\newcommand{\disc}{\mbox{\rm disc}\,}
\newcommand{\eps}{\varepsilon}
\newcommand{\fP}{{\mathfrak P}}
\newcommand{\frf}{\mathfrak f}
\newcommand{\too}{\longmapsto}
\newcommand{\lra}{\longrightarrow}
\newcommand{\Lra}{\Longrightarrow}
\newcommand{\gen}{{\operatorname{gen}}}
\newcommand{\cen}{{\operatorname{cen}}}
\newcommand{\dl}{{\mathcal L}}

\newcommand{\srt}[2]{\!\sqrt[\raisebox{0.2ex}{$\scriptstyle{#1\,\,}$}]{#2}\,}
\newcommand{\im}{\operatorname{im}}
\newcommand{\ba}{\overline{a}}
\newcommand{\bH}{\overline{H}}
\newcommand{\vc}{\mbox{\em v}}

\newtheorem{lem}{Lemma}
\newtheorem{prop}{Proposition}
\newtheorem{thm}{Theorem}

\title{The Development of the Principal Genus Theorem}
\author{Franz Lemmermeyer}
\address{CSU San Marcos, Dept. Mathematics \\
333 S Twin Oaks Valley Rd \\
San Marcos, CA 92096-0001 \\ USA}
\email{franzl@csusm.edu}

\begin{document}
\maketitle

\section*{Introduction}

Genus theory belongs to algebraic number theory and, in
very broad terms, deals with the part of the ideal class 
group of a number field that is `easy to compute'. 
Historically, the importance of genus theory stems from
the fact that it was the essential algebraic ingredient in 
the derivation of the classical reciprocity laws -- from 
Gau\ss{}'s second proof over Kummer's contributions up to 
Takagi's `general' reciprocity law for $p$-th power residues.

The central theorem in genus theory is the principal genus
theorem, which is hard to describe in just one sentence -- 
readers not familiar with genus theory might want to glance
into Section \ref{SG} before reading on. In modern terms, the 
principal genus theorem for abelian extensions $k/\Q$ describes 
the splitting of prime ideals of $k$ in the genus field $k_\gen$ 
of $k$, which by definition is the maximal unramified extension 
of $k$ that is abelian over $\Q$.

In this note we outline the development of the principal
genus theorem from its conception in the context of binary
quadratic forms by Gau\ss{} (with hindsight, traces of 
genus theory can be found in the work of Euler on quadratic
forms and idoneal numbers) to its modern formulation within
the framework of class field theory. It is somewhat remarkable
that, although the theorem itself is classical, the name
`principal ideal theorem' (`Hauptgeschlechtssatz' in German)
was not used in the 19th century, and it seems that it was
coined by Hasse in his Bericht \cite{Has} and adopted immediately
by the abstract algebra group around Noether. It is even more
remarkable that Gau\ss{} doesn't bother to formulate the
principal genus theorem except in passing: after observing
in \cite[\S 247]{Gauss} that duplicated classes (classes of 
forms composed with themselves) lie in the principal genus, 
the converse (namely the principal genus theorem) is stated 
for the first time in \S 261: 
\begin{quote}
si itaque omnes classes generis principalis ex duplicatione 
alicuius classis provenire possunt (quod revera semper locum 
habere in sequentibus demonstrabitur), \ldots\footnote{
if therefore all classes of the principal genus result from the
duplication of some class (and the fact that this is always true
will be proved in the sequel), \ldots} 
\end{quote}
The actual statement of the principal genus theorem is 
somewhat hidden in \cite[\S 286]{Gauss}, where Gau\ss{} formulates
the following
\begin{quote}
Problem. Given a binary form $F = (A, B, C)$ of determinant $D$
belonging to a principal genus: to find a binary form $f$ from
whose duplication we get the form $F$.
\end{quote}
It strikes us as odd that Gau\ss{} didn't formulate this central 
result properly;\footnote{Actually, his formulation is vaguely 
reminiscent of the way in which Euclid presented some of his results.}
yet he knew exactly what he was doing \cite[\S 287]{Gauss}:
\begin{quote}
Since by the solution of the problem of the preceding article
it is clear that any properly primitive (positive) class of
binary forms belonging to the principal genus can be derived
from the duplication of any properly primitive class of the
same determinant, \ldots 
\end{quote}
and he clearly saw the importance of this result:
\begin{quote}
We believe that these theorems are among the most beautiful 
in the theory of binary forms, especially because, despite 
their extreme simplicity, they are so profound that a rigorous 
demonstration requires the help of many other investigations.
\end{quote}

Gau\ss{}'s theory of quadratic forms was generalized in several
completely different directions: 
\begin{enumerate}
\item The theory of $n$-ary quadratic forms over fields, which was 
  cultivated by Hermite, Smith, and Minkowski, and blossomed in the 
  20th century under the hands of Hasse, Witt, Siegel, and 
  others.\footnote{See Jones \cite{Jon}, Lam \cite{Lam}, and O'Meara 
  \cite{OM} for $n$-ary forms, and Buell \cite{Buell} for the binary 
  case; Venkov \cite{Ven} gives a very readable presentation of 
  Gau\ss's results close to the original.} 
\item The arithmetic of algebraic tori encompasses the theory 
  of binary quadratic forms: see Shyr \cite{Shyr1,Shyr2} for 
  a presentation of Gau\ss's theory in this language, and Ono 
  \cite{Ono} for a derivation of the principal genus theorem 
  using results from Shyr's thesis. 
\item The theory of forms of higher degree, in particular cubic 
  forms. We will be content with mentioning only two contributions
  to the algebraic theory of cubic forms: Eisenstein proved several 
  results on cubic forms that nowadays would be presented in the 
  language of cyclic cubic fields (see Hoffman and Morales \cite{HM} 
  for a modern interpretation of composition of cubic forms \`a la Kneser); 
  Manin \cite{Man} studied cubic forms from the viewpoint of obstructions 
  to the local-global principle, and his ideas led to profound 
  insights in modern arithmetic geometry (see Skorobogatov \cite{Skor}). 
\item The theory of quadratic and, later, general algebraic number 
  fields, with Kummer, Dirichlet, Dedekind, and Weber being 
  responsible for the transition from forms to ideal classes in 
  number fields. 
\end{enumerate}

This article is restricted to the genus theory of number fields; 
for a related survey with an emphasis on the quadratic case, but 
sketching generalizations of the genus concept e.g. in group 
theory, see Frei \cite{Frei}. 

\vskip 1cm

\begin{center}
{\huge I. Genus Theory of Quadratic Forms}
\end{center}

\section{Prehistory: Euler, Lagrange and Legendre}

There are hardly any traces of genus theory in the mathematical 
literature prior to Gau\ss{}'s Disquisitiones. What can be found, in 
particular in Euler's work, are results and conjectures that later
were explained by genus theory.

One such conjecture was developed between Goldbach and Euler:
on March 12, 1753, Goldbach wrote to Euler \cite[Letter 166]{EG} 
that if $p$ is a prime $\equiv 1 \bmod 4d$, then $p$ can be 
represented as $p = da^2 + b^2$. Euler replies on March 23/April 3 
\cite[Letter 167]{EG}:
\begin{quote}
Ich habe auch eben diesen Satz schon l\"angst bemerket und bin von
der Wahrheit desselben so \"uberzeugt, als wann ich davon eine 
Demonstration h\"atte.\footnote{I have known this very theorem for
quite a long time, and I am just as convinced of its truth as if 
I had its demonstration.}
\end{quote} 
He then gives the examples
\begin{align*}
 p & = 4 \cdot 1m + 1 & \Lra & & p & =\ aa + bb \\
 p & = 4 \cdot 2m + 1 & \Lra & & p & = 2aa + bb \\
 p & = 4 \cdot 3m + 1 & \Lra & & p & = 3aa + bb \\
 p & = 4 \cdot 5m + 1 & \Lra & & p & = 5aa + bb \quad etc.
\end{align*}
and remarks that he can prove the first claim, but not the
rest.\footnote{Later he found a proof for the case $p = 3a^2 + b^2$;
the other two cases mentioned here were first proved by Lagrange.}
Euler then goes on to observe that the conjecture is only true in 
general when $a$ and $b$ are allowed to be rational numbers, and 
gives the example $89 = 4 \cdot 22 + 1$, which can be written as 
$89 = 11(\frac{5}{2})^2 + (\frac{9}{2})^2$ but not in the form
$11a^2 + b^2$ with integers $a, b$. Thus, he says, the theorem
has to be formulated like this:

\medskip\noindent{\bf Conjecture 1.}
{\em Si $4n+1$ sit numerus primus, et $d$ divisor ipsius $n$, tum 
iste numerus $4n+1$ certo in hac forma $daa+bb$ continentur, si non 
in integris, saltem in fractis.}\footnote{If $4n+1$ is a prime 
number, and $d$ a divisor of this $n$, then that number $4n+1$ is
certainly contained in the form $daa+bb$, if not in integers,
then in fractions.}
\medskip

Euler also studied the prime divisors of a given binary quadratic
form $x^2 + ny^2$ (in the following, we will always talk about 
`proper' divisors of quadratic forms, that is, we assume that 
$p \mid x^2 + ny^2$ with $\gcd(x,y) = 1$), and observed that 
those not dividing $4n$ are contained in half of the possible 
residue classes modulo $4n$ coprime to $4n$: for example, the 
prime divisors of $x^2 + 5y^2$ are contained in the residue 
classes $1, 3, 7, 9 \bmod 20$; no prime congruent to 
$11, 13, 17, 19 \bmod 20$ divides $x^2 + 5y^2$ without dividing
$x$ and $y$.\footnote{Euler \cite[p. 210]{EulI} expressed this 
by saying that primes (except $p = 2, 5$) dividing $x^2 + 5y^2$ 
have the form $10i \pm 1$, $10i \pm 3$, where the plus sign holds 
when $i$ is even, and the minus sign when $i$ is odd.} Now, as 
Euler knew (and he used this in his `solution' of the cubic Fermat 
equation), odd primes dividing $x^2+ny^2$ can be represented by 
the same quadratic form if $n=3$, and he also knew that this 
property failed for $n = 5$. He then saw that the primes
$p \equiv 1, 9 \bmod 20$ could be represented\footnote{At this
stage, he had already studied Lagrange's theory of reduction of 
binary quadratic forms.} as $p = x^2 + 5y^2$ with $x, y \in \N$, 
whereas $p \equiv 3, 7 \bmod 20$ could be written as 
$2x^2 + 2xy + 3y^2$ with $x, y \in \Z$.  His first guess was that 
this would generalize as follows: the residue classes containing 
prime divisors of $x^2 + ny^2$ could be associated uniquely with a 
reduced quadratic form of the same discriminant as  $x^2 + ny^2$.
For example, the reduced forms associated to $F = x^2 + 30y^2$ are
the forms $D$ satisfying $D = F$, $2D = F$, $3D = F$ and $5D = F$,
where $2D = F$ refers to $D = 2r^2 + 15 s^2$, $3D = F$ to 
$3r^2 + 10s^2$, and $5D = F$ to $D = 5r^2 + 6s^2$. Each of these
forms has different classes of divisors. 

As Euler \cite[p. 192]{EulI} found out, however, $n=39$ provides 
a counterexample: he comes up with three classes of forms
\begin{quote}
Hinc igitur patet omnino dari tria genera divisorum:\footnote{Thus 
from here it becomes clear that there are altogether three kinds of
divisors:}
$$ 1)\ D = F, \qquad 2)\ 3D = F, \qquad [3)]\ 5D = F. $$
\end{quote}
These three kinds of divisors are $D = F = r^2 + 39s^2$; 
$D = 3r^2 + 13s^2$ (note that $3D = (3r)^2 + 39s^2$, which 
explains Euler's notation $3D = F$); and $D = 5r^2 + 2rs + 8s^2$. 
He then observes that the divisors of the first and the second 
class share the same residue classes modulo $156$; the prime 
$61 = 3 \cdot 4^2 + 13 \cdot 1^2$ belonging to the second 
class can be represented rationally by the first form since 
$61 = (\frac{25}{4})^2 + 39 (\frac{3}{4})^2$.

\subsection*{Euler and Lagrange}
One of the results in which Euler comes close to genus theory
is related to a conjecture of Euler that was shown to be false 
by Lagrange; it appears in \cite{EulR}. In his comments on 
Euler's algebra, Lagrange \cite[p. 156--157]{Lagr} writes
\begin{quote}
Euler, dans un excellent M\'emoire imprim\'e dans le
tome IX des {\em nouveaux Commentaires de P\'etersbourg},
trouve par induction cette r\`egle, pour juger la
r\'esolubilit\'e de toute equation de la forme
$$ x^2 - Ay^2 = B,$$
lorsque $B$ est un nombre premier; c'est que l'\'equation
doit \^etre possible toutes les fois que $B$ sera de la
forme $4An+r^2$, ou $4An+r^2-A$;\footnote{M. Euler, in an
excellent Memoir printed in vol. IX of the {\em New
Commentaries of Petersburg}, finds by induction this rule
for determining the solvability of every equation of the
form $x^2 - Ay^2 = B$, where $B$ is a prime number:
the equation must be possible whenever $B$ has the form
$4An+r^2$, or $4An + r^2 - A$;}
\end{quote}

For example, $-11 = 4 \cdot 3 \cdot (-1) + 1^2$, and
$-11 = 1^2 - 3 \cdot 2^2$. Similarly,
$-2 =  4 \cdot 3 \cdot (-2) + 5^2 - 3$ and $-2 = 1^2 - 3 \cdot 1^2$.
Euler's main motivation for this conjecture were numerical data,
but he also had a proof that $p = x^2 - ay^2$
implies $p = 4an+r^2$ or $p = 4an + r^2 - a$. In fact, he 
writes $x = 2at+r$, $y = 2q+s$, and finds that $p = x^2 - ay^2 
= 4am + r^2 - as^2$ for some $m \in \Z$. If $s$ is even, then
$-as^2$ has the form $4am'$, and if $s$ is odd, we find 
$-as^2 = -4am''-a$. This proves the claim. 

As Lagrange pointed out, however, Euler's conjecture is not 
correct, and he came up with the following counterexample: 
the equation $x^2 - 79y^2 = 101$ is not solvable in integers, 
although $101 = 4An+r^2-A$ with $A = 79$, $n = -4$ and $r = 38$.

Whether Euler ever heard about Lagrange's counterexample 
is not clear; based on Euler's experience in such matters 
it is not unreasonable to suspect that he probably would have 
reacted in the same way as in the other cases above, namely by
replacing the representation in integers by representation in 
rational numbers. This would have led to the following

\medskip\noindent{\bf Conjecture 2.}
{\em If $p \nmid 4a$ is a prime of the form $4an+r^2$ or 
$4an+r^2-a$, then $p = x^2 - ay^2$ for rational numbers $x, y$.}
\medskip

As we shall see, this conjecture is equivalent to Gau\ss's
principal genus theorem.


\subsection*{Euler and Genus Theory}
Antropov \cite{Ant1,Ant2,Ant3} has tried to make a case for 
the claim that the concept of genera is due to Euler. In 
\cite{Ant1}, he writes
\begin{quote}
It seems to have gone unnoticed that Euler partitioned the set 
of the integral binary quadratic forms with a given discriminant
into classes, which he called ``genera''.
\end{quote}
Euler did use the word `genus', as we have seen, for distinguishing
between certain kinds of quadratic forms; however, Euler uses both 
genus (``ad genus tertium pertinentis'') and class (``ad tertiam 
classem pertinet'') when talking about what Antropov perceives as 
a `genus' of quadratic forms. 

\medskip 

Fueter, in his preface \cite[p. xiii]{Fue}, praises 
Euler's contributions as follows:
\begin{quote}
Man kann sagen, da\ss{} eigentlich alle Bausteine der Zahlentheorie
dieses Gebietes in den {\sc Eulerschen} Abhandlungen schon 
bereitgestellt wurden, und es nur der Meisterhand eines {\sc Gauss}
bedurfte, um sie zu dem Geb\"aude der Theorie der quadratischen 
Formen zusammenzuf\"ugen.\footnote{One might argue that essentially 
all the building blocks of the number theory of this area had been
provided by {\sc Euler}, and that it only took the masterly hand of 
{\sc Gauss} to join them to the edifice of the theory of binary
quadratic forms.}
\end{quote} 
A more realistic view was later offered by Weil, who -- referring
to Euler's papers on idoneal numbers -- wrote (\cite[p. 224]{Weil})
\begin{quote}
They are [\ldots] ill coordinated with one another, and some of
the formulations and proofs in them are [\ldots] confused and
defective \ldots
\end{quote}
One must not forget, however, that Euler only had isolated 
results on (divisors of) numbers represented by quadratic forms,
and that it may be unfair to judge his work in the light of the
insight provided by Gau\ss{}, who subsumed Euler's (and Lagrange's) 
results into just a few theorems (reciprocity, class group,
principal genus theorem) within his theory of quadratic forms.

\medskip

\subsection*{Legendre}
As for Legendre's contribution, Dirichlet \cite[p. 424]{Dir} writes
\begin{quote}
Les formes diff\'erentes qui correspondent au d\'eterminant 
quelconque $D$, sont divis\'ees par M. {\sc Gauss} en genres, 
qui sont analogues \`a ce que {\sc Legendre} appelle groupes 
des diviseurs quadratiques.\footnote{The different forms that
correspond to some determinant $D$ are partitioned by Mr.
{\sc Gauss} into genera, which are analogous to what 
{\sc Legendre} calls groups of quadratic divisors.}
\end{quote}

Legendre's `diviseurs quadratiques' of a binary quadratic forms
are the reduced nonequivalent quadratic forms of the same 
`determinant'; to each of these classes he associates `diviseurs
lineaires', namely the linear forms $ax+b$ with the property that
the primes represented by the `diviseur quadratique' are contained
in the arithmetic progression $ax+b$ (see \cite[Art. 212]{Leg} for
the $8$ diviseurs quadratiques of $x^2 + 41y^2$ and the $6$ diviseurs
lineaires corresponding to each of them).

\section{Gau\ss{}}\label{SG}

Let us start by briefly recalling Gau\ss{}'s definitions. In Section
V of his Disquisitiones Arithmetica, he studies binary quadratic 
forms $F(x,y) = ax^2 + 2bxy + cy^2$ that are occasionally denoted 
by $(a,b,c)$; the {\em determinant} of $F$ is $D = b^2 - ac$. An 
integer $n$ is said to be {\em represented} by $F$ is there exist 
integers $x, y$ such that $n = F(x,y)$. A form $(a,b,c)$ is 
{\em ambiguous} if $a \mid 2b$, and {\em primitive} if 
$\gcd(a,b,c) = 1$. 

The following theorem proved in \S 229 is the basis for the 
definition of the genus of a binary quadratic form:
\begin{quote}
Let $F$ be a primitive binary quadratic form with determinant $D$, 
and let $p \mid D$ be prime. Then the numbers not divisible by $p$ 
that can be represented by $F$ agree in that they are either all
quadratic residues of $p$, or they are all nonresidues.
\end{quote}
For $p = 2$ the claim is correct but trivial. If $4 \mid D$, however,
then the numbers represented by $f$ are all $\equiv 1 \bmod 4$,
or all $\equiv 3 \bmod 4$. Similarly, if $8 \mid D$, the numbers
lie in exactly one of the four residue classes $1, 3, 5$ or $7 \bmod 8$.

Gau\ss{} adds without proof the remark that there is no such pattern 
for odd primes {\em not} dividing the discriminant:

\medskip\noindent{\bf Observation.}\label{GO}
{\em If it were necessary for our purposes we could easily show 
that numbers representable by the form $F$ have no such 
fixed relationship to a prime number that does not divide 
$D$} \footnote{\"Ubrigens w\"urden wir, wenn es zum gegenw\"artigen 
Zwecke notwendig w\"are, leicht beweisen k\"onnen, da\ss{} die 
durch $F$ darstellbaren Zahlen zu keiner in $D$ nicht aufgehenden 
Primzahl in einer derartigen festen Beziehung stehen.}  

\medskip

except for the residue classes modulo $4$ and $8$ of representable
odd numbers in case $D$ is odd:
\begin{enumerate} 
\item[I.] If $D \equiv 3 \bmod 4$, then the odd $n$ that can be 
          represented by $F$ are all $\equiv 1 \bmod 4$ or all 
          $\equiv 3 \bmod 4$.
\item[II.] If $D \equiv 2 \bmod 8$, then the odd $n$ that can be 
          represented by $F$ are all  $\equiv \pm 1 \bmod 8$ or
          all $\equiv \pm 3 \bmod 8$.
\item[III.] If $D \equiv 6 \bmod 8$, then the odd $n$ that can be 
          represented by $F$ are all $\equiv 1, 3 \bmod 8$ or all
          $\equiv 5, 7 \bmod 8$.
\end{enumerate}


Gau\ss{} uses this observation to define characters of primitive 
quadratic forms (\S 230); for example, to the quadratic form
$(7, 0, 23)$ he attaches the {\em total character} $1,4; R7; N23$ 
because the integers represented by $7x^2 + 23y^2$ are 
$\equiv 1 \bmod 4$, quadratic residues modul $7$, and 
quadratic nonresidues modulo $23$. Gau\ss{} observes 
that if $(a,b,c)$ is a primitive quadratic form, then 
$p \mid b^2-ac$ implies $p \nmid \gcd(a,c)$, so the character of 
primitive forms can be determined from the integers $a$ and $c$, 
which of course are both represented by $(a,b,c)$. Finally he
remarks that forms in the same class have the same total character,
which allows him to consider them as characters of the classes.
   
Now Gau\ss{} collects classes of forms of given determinant into 
genera: a genus is simply the set of all classes with the same
total character. The principal genus is the genus containing
the principal class (the class containing the principal form 
$(1, 0, -D)$). 

Next on Gau\ss{}'s agenda (\S 234 -- \S 256) is the definition of 
the composition of forms, orders, genera,\footnote{This terminology 
is apparently taken from biology. C.~Linne classified the living 
organisms into kingdoms (plants and animals), classes, orders, genera, 
and species. Kummer used the German expression `Gattung' for Gau\ss's 
`genus', but in the long run the translation `Geschlecht' prevailed.} 
and classes of forms. The next three articles (\S 257 -- \S 259) are
devoted to the determination of the number of ambiguous classes.

In \S 261 Gau\ss{} proves the {\em first inequality} of genus theory:
at least half of all possible total characters do not occur. This
is a consequence of the `ambiguous class number formula'. \S 262 is 
reserved for a demonstration that the first inequality implies the 
quadratic reciprocity law. 

In \S 266, Gau\ss{} begins a long excursion into the theory of 
ternary quadratic forms $Ax^2 + 2Bxy + Cy^2 + 2Dxz + 2Eyz + Fz^2$.
After discussing the reduction of ternary quadratic forms, 
he remarks that ternary forms represent both integers (by 
substituting integers for $x$, $y$, $z$) and quadratic forms
(by putting e.g. $x = at+bu$, $y = ct+du$, $z = et+fu$ for
variables $t$ and $u$). 

After having studied the representations of binary quadratic forms
by ternary forms, Gau\ss{} returns to binary quadratic forms in 
\S 286 and proves the principal genus theorem: every form $F$ in 
the principal genus is equivalent to $2f$ for some form $f$ of the 
same determinant as $F$. This immediately implies the {\em second 
inequality} of genus theory in \S 287: at least half of all possible 
total characters do in fact occur. Finally, in \S 303, Gau\ss{} 
characterizes Euler's idoneal numbers using genus theory. 

\subsection*{Ambiguous Forms and Classes.}
In his Marburg lectures, Hasse \cite[p. 158]{HasM} writes
\begin{quote}
Die in diesem Zusammenhang nicht sehr gl\"uckliche Bezeichnung
``ambig'' stammt von Gau\ss.\footnote{The term ``ambiguous'', whose 
usage in this connection is somewhat unfortunate, is due to Gau\ss{}.}
\end{quote}
Gau\ss{}, however, wrote in Latin and used `forma anceps' to 
denote an ambiguous quadratic form; apparently Hasse confused 
the disquisitiones with Maser's translation, where `anceps' 
is translated as `ambig'.\footnote{Clarke used `ambiguous' 
in his English translation, and I.~Adamson used `ambig' 
in his English translation of Hilbert's Zahlbericht.}

The actual story of the word ambiguous is given by
Dedekind \cite[\S 58, p. 139]{DD}:
\begin{quote}
In his lectures Dirichlet always used the word {\em forma anceps},
which I have kept when I prepared the first edition (1863); in the
second and third edition (1871, 1879), [\ldots] I called them
{\em ambiguous forms} following Kummer, who used this notation 
in a related field; 
\end{quote} 
apparently there were complaints about the word, and in the
fourth edition he replaces `ambiguous' by `twosided' (`zweiseitig'
in German). 

As a matter of fact, the expression `ambiguous' was not at all
Kummer's invention: it was used in the form `classe ambigu\"e' by 
Poullet Delisle in his French translation of the disquisitiones, 
which appeared in 1807. By the time Kummer started studying number
theory, the Latin edition of the disquisitiones must have been
next to impossible to get; we know that Eisenstein's copy of 
the disquisitiones was in French (see Weil \cite{WE}), and it
seems reasonable to assume that Kummer studied the same edition.
Nowadays, the French use the word `ambige'; apparently it was 
Chevalley \cite{Chev1} who dropped the `u'.\footnote{I owe this 
remark to J.-F. Jaulent.}

\subsection*{Ternary Forms}
Gau\ss{} used his theory of ternary quadratic forms to prove the 
principal genus theorem, and derived Legendre's theorem (as well as 
the $3$-squares theorem\footnote{Every positive integer not of the 
form $4^a(8b+7)$ can be written as a sum of three squares.}) from the 
same source. Arndt \cite{Arndt} and later Dedekind \cite[\S 158]{DD}
and Mansion \cite{Mans} realized that Legendre's theorem is sufficient 
for proving the principal genus theorem, which simplified the theory 
considerably (see \cite[Chap. 2]{L1}). 
Unfortunately, however, Legendre's theorem does not seem to suffice 
for deriving the $3$-squares theorem, but Deuring \cite[VII, \S 9]{Deu} 
(see also Weil \cite[III, App. II; p. 292--294]{Weil}) sketched a 
very beautiful proof using the theory of quaternion algebras. 
Venkov (1927; see \cite{Ven}) used Gau\ss's theory of ternary 
quadratic forms to give an arithmetic proof of Dirichlet's class 
number formula for negative discriminants $-m$ in which $m$ is the 
sum of three squares.  Shanks \cite{ShaG} used binary quadratic 
forms to develop his extremely clever factorization algorithm 
SQUFOF,\footnote{SQUare FOrm Factorization.} and Gau\ss's theory 
of ternary quadratic forms \cite{ShaT} for coming up with an 
algorithm for computing the $2$-class group of complex quadratic 
number fields.

\section{Dirichlet-Dedekind}

According to a well known story (see Reichardt \cite[p. 14]{Rei}), 
Dirichlet never put Gau\ss{}'s disquisitiones on the bookshelf
but kept the copy on his desk and took it with him on journeys.
Dirichlet's constant occupation with the disquisitiones provided
him with the insight that allowed him to streamline and simplify
Gau\ss's exposition (occasionally by restricting himself to a 
special case), thereby making the disquisitiones accessible 
to a much wider audience. 

In \cite{Dir}, Dirichlet replaces Gau\ss's notation $aRp$ by
$(a/p) = +1$, thus giving Gau\ss's characters the now familiar 
look. His main contribution in \cite{Dir} was definitely the 
proof of the `second inequality' of genus theory using analytic 
methods.\footnote{See Zagier \cite{Zag} for a modern exposition.} 

Dirichlet presented the theory of binary quadratic forms in his 
lectures; his results on genus theory were added by Dedekind in 
the supplements IV (analytic proof) and X (arithmetic proof using 
Legendre's theorem). In \S 122, he defines an integer 
$$\lambda = \# \{\text{odd primes dividing $D$}\}\ + \begin{cases}
     0 & \text{if}\ D \equiv 1 \bmod 4 \\
     2 & \text{if}\ D \equiv 0 \bmod 8 \\
     1 & \text{otherwise}, \end{cases} $$ 
and in \S 123 he proves the first inequality of genus theory:
$$ g \le 2^{\lambda-1}. $$
In \S 125, he gives Dirichlet's analytic proof of the existence
of these genera, the second inequality of genus theory:
\begin{quote} 
Die Anzahl der wirklich existierenden Geschlechter ist gleich 
$2^{\lambda-1}$, und alle diese Geschlechter enthalten gleich 
viele Formenklassen.\footnote{The number of existing genera is
$2^{\lambda-1}$, and all these genera contain equally many 
classes of forms.}
\end{quote}
He also remarks that the second inequality would follow immediately
from Dirichlet's theorem on the infinitude of primes in arithmetic
progressions.

Dedekind returns to genus theory of binary quadratic forms in 
his supplement X: \S 153 gives the first inequality, \S 154 the 
quadratic reciprocity law, and in \S 155 he observes that the
second inequality of genus theory (the existence of half of all
the possible genera) is essentially identical with the 
principal genus theorem: 
\begin{quote}
Every class of the principal genus arises from duplication.
\end{quote}
He then remarks
\begin{quote}
Wir k\"onnen hier unm\"oglich darauf eingehen, den Beweis 
mit\-zu\-thei\-len, welchen Gauss auf die Theorie der 
tern\"aren quadratischen Formen gest\"utzt hat; da dieses 
tiefe Theorem aber den sch\"onsten Abschluss der Lehre von 
der Composition bildet, so k\"onnen wir es uns nicht versagen, 
dasselbe auch ohne H\"ulfe der Dirichlet'schen Principien auf 
einem zweiten Wege abzuleiten, der zugleich die Grundlage f\"ur
andere wichtige Untersuchungen bildet.\footnote{It is impossible
for us to communicate the proof, which Gauss has based on the
theory of ternary quadratic forms; but since this deep theorem
is the most beautiful conclusion of the theory of composition, we 
cannot help but derive this result, without the use of Dirichlet's 
principles,  in a second way, which will also form the basis for
other important investigations.}
\end{quote}
His proof begins by showing that the following statement is
equivalent to the principal genus theorem:
\begin{quote}
If $(A,B,C)$ is a form in the principal genus of determinant $D$,
then the equation
$$ Az^2 + 2Bzy + Cy^2 = x^2 $$
has solutions in integers $z, y, x$ such that $x$ is coprime to $2D$.
\footnote{Ist $(A,B,C)$ eine Form des Hauptgeschlechtes der Determinante 
$D$, so ist die Gleichung 
\begin{equation}\label{EDd} Az^2 + 2Bzy + Cy^2 = x^2 \end{equation}
stets l\"osbar in ganzen Zahlen $z, y, x$, deren letzte relative
Primzahl zu $2D$ ist.}
\end{quote}
In  \S 158 Dedekind gives a proof of the principal genus theorem 
based on Legendre's theorem (which, as he observes in a footnote 
in \S 158, belongs to the theory of ternary quadratic forms) and 
refers to Arndt \cite{Arndt} for a first proof of this kind.

\medskip
\noindent{\bf Remark.}
The genus theory of Dirichlet and Dedekind is a genus theory of
binary quadratic forms. Although Dedekind introduced ideals and
maximal orders in number fields, he did not translate genus theory
into his new language.

Dirichlet's analytic methods were used to prove the principal 
genus theorem by Kronecker \cite{Kro} and de S\'eguier 
\cite[p. 135--153; 333-334]{dS}. 
A proof based on Dedekind's criterion involving (\ref{EDd}) was 
given by P\'epin \cite[p. 45]{Pep}.

Mertens \cite{Mer} gave a new proof of the principal genus theorem
built on Legendre's theorem. De la Vall\'ee Poussin \cite{VP} 
and Mertens \cite{Mer2} found proofs based only on binary quadratic 
forms. See also Speiser \cite{Sp}.

Heine \cite{Hei} studied quadratic forms over function fields
of one variable; see also Bae \& Koo \cite{BaKo} and 
Hellegouarche \cite{Hel}. The `genre principal' introduced by 
Serret \cite{Ser} in his investigations of irreducible polynomials 
of degree $p^\mu$ over $\F_p[X]$ seems unrelated to the principal 
genus of Gau\ss.

Dirichlet \cite{Diri}, Smith \cite{Smith}, and Minnigerode \cite{Min}
investigated binary quadratic forms with coefficients in $\Z[i]$. 
Speiser \cite{Spei} developed genus theory for binary quadratic
forms with coefficients from the ring of integers of an arbitrary
number field.
\vskip 1cm

\begin{center}
{\huge II. Genus Theory of Quadratic Number Fields}
\end{center}

\section{Hilbert} 

Before Hilbert published his report on algebraic numbers, he worked on 
the arithmetic of quadratic extensions of $\Q(i)$ (Dirichlet fields) 
with the intention of 
\begin{quote}
extending the theory of Dirichlet's biquadratic number field in a 
purely arithmetic way to the same level that the theory of quadratic 
number fields has had since {\sc Gauss},
\end{quote}
and the main tool for achieving this goal was, according to Hilbert, 
the notion of genera of ideal classes. \footnote{The complete 
quotation from \cite{Hi} reads: Die vorliegende Abhandlung hat das 
Ziel, die Theorie des Dirichletschen biquadratischen Zahlk\"orpers auf 
rein arithmetischem Weg bis zu demjenigen Standpunkt zu f\"ordern, 
auf welchem sich die Theorie der quadratischen K\"orper bereits seit 
{\sc Gauss} befindet. Es ist hierzu vor allem die Einf\"uhrung des 
Geschlechtsbegriffs sowie eine Untersuchung derjenigen Einteilung 
aller Idealklassen notwendig, welche sich auf den Geschlechtsbegriff 
gr\"undet.}

Let $\Z[i]$ denote the ring of Gaussian integers, and let 
$\delta \in \Z[i]$ be a squarefree nonsquare. Hilbert considers
the quadratic extension $K = \Q(\sqrt{\delta}\,)$ of $k = \Q(i)$,
computes integral bases, and determines the decomposition of
primes. 

For the definition of the genus Hilbert introduces the prototype of
his norm residue symbol. For $\sigma \in k$ and $\lambda$ a 
prime divisor $\ne (1+i)$ of the discriminant of $K/k$, Hilbert
writes $\sigma = \alpha\nu$ as a product of a relative norm $\nu$
and some $\alpha \in \Z[i]$ not divisible by $\lambda$, and puts
$$ \Big[ \frac{\sigma}{\lambda:\delta}\Big] = 
   \Big[ \frac{\alpha}{\lambda}\Big], $$
where $[\,\cdot\,/\,\cdot\,]$ is the  quadratic residue symbol 
in $\Z[i]$. The definition for $\lambda = 1+i$ is slightly more
involved.

Then Hilbert defines the character system of an ideal $\fa$ in $\OO_K$
as the system of signs
$$ \Big[ \frac{\sigma}{\lambda_1:\delta}\Big] , \ldots, 
   \Big[ \frac{\sigma}{\lambda_s:\delta}\Big], $$
where $\lambda_1$, \ldots, $\lambda_s$ denote the ramified primes.
The character system of ideals only depends on their ideal class, 
and classes with the same character system are then said to be in 
the same genus. The principal genus is the set of ideal classes
whose character system is trivial. The principal genus theorem
is formulated in \cite[\S 4]{Hi}:
\begin{quote}
Each ideal class in the principal genus is the square of some
ideal class.\footnote{Eine jede Idealklasse des Hauptgeschlechtes 
ist gleich dem Quadrat einer Idealklasse.}
\end{quote}
Hilbert then determines the number of genera, derives the quadratic
reciprocity law, and finally gives an arithmetic proof of the
class number formula for $\Q(i,\sqrt{m}\,)$ and $m \in \Z$.
He apparently has not yet realized that his symbols 
$\big[\frac{\sigma}{\lambda:\delta}\big]$ are `norm residue' symbols,
or that the quadratic reciprocity law can be expressed by a 
product formula.

He takes these steps in the third section of his Zahlbericht,
which deals with the theory of quadratic number fields. He
calls an integer $n$ a norm residue\footnote{I will adapt the 
following convention: an element is a norm residue {\em modulo} 
$\fa$ if it is congruent to a norm modulo $\fa$, and a norm residue
at $\fp$ if it is congruent to norms modulo any power $\fp^k$.} 
at $p$ in $\Q(\sqrt{m}\,)$ if $m$ is a square or if for all 
$k \ge 1$ there exist integers $x, y \in \Z$ such that 
$n \equiv x^2 - my^2 \bmod p^k$. Then he defines the norm 
residue symbol by
$$ \Big(\frac{n\,,\,m}{p}\Big) = \begin{cases}
   +1 & \text{if $m$ is a norm residue at $p$ in $\Q(\sqrt{m}\,)$} \\
   -1 & \text{otherwise}.   
   \end{cases} $$
Hilbert uses the norm residue symbol to define characters on 
ideal classes and defines the principal genus to consist of 
those ideal classes with trivial character system. In \S 68,
he employs ambiguous ideals and his Satz 90 to prove that 
quadratic number fields with exactly one ramified prime have 
odd class number, and then deduces the quadratic reciprocity 
law in \S 69. In \S 72 he proves the principal genus theorem:
\begin{quote}
In einem quadratischen K\"orper ist jede Klasse des Hauptgeschlechts
stets gleich dem Quadrat einer Klasse ({\sc Gauss}).\footnote{In a 
quadratic number field, each class of the principal genus is the
square of a class ({\sc Gauss}).}
\end{quote}
The proof uses a reduction technique reminiscent of Lagrange; 
the solvability of the `norm equation' $n = x^2 - my^2$ for 
$x, y \in \Q$ is equivalent to the fact that the ternary
quadratic form $x^2 - my^2 - nz^2$ nontrivially represents 
$0$ in integers, and Hilbert explicitly refers to Lagrange when
he states the following special case of `Hasse's norm theorem':
\begin{quote}
{\sc Satz} 102. Wenn $n, m$ zwei ganze rationale Zahlen bedeuten,
von denen $m$ keine Quadratzahl ist, und die f\"ur jede beliebige
Primzahl $w$ die Bedingung
$$ \Big(\frac{n\,,\,m}{w}\Big) = +1 $$
erf\"ullen, so ist die Zahl $n$ stets gleich der Norm 
einer ganzen oder gebrochenen Zahl $\alpha$ des K\"orpers 
$k(\sqrt{m}\,)$.\footnote{{\sc Theorem} 102. If $n, m$ denote
two rational integers, where $m$ is a nonsquare, and if for
any prime $w$ the condition  
$$ \Big(\frac{n\,,\,m}{w}\Big) = +1 $$
is satisfied, then $n$ is the norm of a (not necessarily integral)
number $\alpha$ of the field $k(\sqrt{m}\,)$.}
\end{quote}
(Note that $k(\sqrt{m}\,)$ denotes the quadratic number field $k$
one gets by adjoining $\sqrt{m}$ to the field of rational numbers.)

The ambiguous class number formula is proved afterwards, and finally
Hilbert gives a second proof of the principal genus theorem using
Dirichlet's analytic techniques.

\medskip
\noindent{\bf Remark.}
With Hilbert, the transition from Gau\ss's genus theory of binary
quadratic forms to the corresponding theory of quadratic extensions
is complete. Distinctive features of Hilbert's presentation are
\begin{enumerate} 
\item the central role of the ambiguous class number formula; 
       the cohomological kernel of these results was recognized 
       only much later;
\item the introduction of norm residue symbols and the formulation
       of the reciprocity law as a product formula.
\end{enumerate}
Although Hilbert saw that the norm residue symbol for the `infinite 
prime' of $\Q$ (he wrote it as $(\frac{n,m}{-1})$; see \cite[\S 70]{Hil}) 
would simplify the presentation, he chose not to use it. These symbols
became necessary when he replaced $\Q$ by arbitrary base fields $k$
in his article \cite{HilA} on class field theory.

\section{Weber}

In the third volume of his algebra \cite[\S 108]{WA3}, Weber 
gives an account of genus theory that shows Hilbert's influence: 
while Weber does not include the theory of the quadratic Hilbert 
symbol, he realizes the importance of the concept of norm residues.

He fixes a modulus $m \in \N$, considers a natural number $S$
divisible by $m$, and forms\footnote{Presumably this is 
influenced by Hensel.} the multiplicative group $Z$ of
rational numbers coprime to $S$, that is, the set of all 
$\frac{a}{b}$ with $a, b \in \Z$ and $\gcd(a,S) = \gcd(b,S) = 1$. 
The kernel of the natural map $Z \lra (\Z/m\Z)^\times$ is the 
group $M$ of all elements of $Z$ that are congruent to $1 \bmod m$,
and Weber observes that $(Z:M) = \phi(m)$. 

Now let $\OO$ denote an order of a quadratic number field $k$
(Weber writes $Q$ instead of $\OO$) such that the prime factors 
of the conductor of $\OO$ divide $S$; in particular, the discriminant 
$\Delta$ of $\OO$ is only divisible by primes dividing $S$. The set 
of integers $a \in \Z$ for which there is an $\omega \in \OO$ with 
$N \omega \equiv a \bmod m$ form a subgroup $A$ of $M$, the group 
of norm residues modulo $m$ of $\OO$.

In order to simplify the presentation, let $A\{m\}$ denote the
group of norm residues modulo $m$. Weber \cite[\S 107]{WA3} observes 
that if $m = m_1m_2$ with $\gcd(m_1,m_2) = 1$, then
$(Z:A\{m\}) = (Z:A\{m_1\})(Z:A\{m_2\})$, thereby reducing the 
computation of the index $(Z:A\{m\})$ to the case of prime powers $m$.

In \cite[\S 108]{WA3}, he proves
$$ (Z:A\{p^t\}) = \begin{cases}
     1 & \text{if}\ p \nmid \Delta \\
     2 & \text{if}\ p \mid  \Delta \end{cases} $$
if $p$ is an odd prime, and
$$ (Z:A\{2^t\}) = \begin{cases}
     1 & \text{if}\ \Delta \equiv 1 \bmod 4, \ \Delta \equiv 4, 20 \bmod 32, \\
     2 & \text{if}\ \Delta \equiv 8, 12, 16, 24, 28 \bmod 32, \\
     4 & \text{if}\ \Delta \equiv 0 \bmod 32. \end{cases} $$

If $r$ is norm residue modulo $m$ for any modulus $m$ coprime to $r$,
then $r$ is called an absolute norm residue, and the set of all such
$r \in Z$ forms a group $R$ (\cite[\S 109]{WA3}). As a consequence of 
his index computations above, Weber records
$$ (Z:R) = 2^\lambda, $$
where $\lambda$ is the number of discriminant divisors of $\Delta$.
A divisor $\delta$ of $\Delta$ is called a discriminant divisor 
if both $\delta$ and $\Delta/\delta$ are discriminants.

Weber \cite[\S 109]{WA3} defines the genus of an ideal as the set
of all ideals $\fa$ coprime to $\Delta$ whose norms $N\fa$ are in 
the same coset of $Z/R$, and observes that equivalent ideals have 
the same genus. The principal genus is the group of all ideals 
coprime to $\Delta$ such that $N\fa \in R$. He shows that the
existence of primes that are quadratic nonresidues modulo $\Delta$ 
implies that the number $g$ of genera satisfies the inequality
\begin{equation}\label{EWg}  g  \le \frac12 (Z:R), \end{equation}
and that the existence of such primes is equivalent to the 
quadratic reciprocity law (once more we can observe the close
connection between quadratic reciprocity and the first inequality).
The fact that we have equality in (\ref{EWg}) is proved in
\cite[\S 113]{WA3} using Dirichlet's analytic methods.

\medskip
\noindent{\bf Remark.}
Weber's presentation of the genus theory of orders in quadratic
number fields was taken up (to some degree) in Hecke's textbook. 
Observe that the `localization' of the necessary index calculations 
is much more visible in Weber's than in Hilbert's treatment. 
Moreover, these index formulas are intimately related to Gau\ss's
observation on p. \pageref{GO}.

\section{Hecke}

Hecke's `Vorlesungen' \cite{Hecke} contain a masterful exposition 
of algebraic number theory, including the genus theory of (the
maximal orders of) quadratic fields. During the reformulation
of class field theory in the 1930s, genus theory was thrust into 
the background as local methods gradually replaced genus theory 
in the foundation of class field theory. 

Hecke's presentation of genus theory in quadratic fields $k$
with discriminant $d$ combines known with new features:
\begin{itemize}
\item Hecke uses class groups in the strict sense; already Hilbert
      \cite[\S 83--84]{Hil} had seen that this simplifies the 
      exposition of genus theory because some of the statements
      ``can be expressed in a simpler way by using the new 
      notions''.\footnote{\ldots und einige [dieser Tatsachen] 
     erhalten bei Verwendung der neuen Begriffe sogar noch einen 
     einfacheren Ausdruck.}
\item Hecke uses Weber's index computation for norm residues,
      but restricts his attention right from the start to
      norm residues modulo $d$. 
\item Hecke gives a new and very simple definition of genera: two
      ideals $\fa$ and $\fb$ coprime to $d$ are said to belong to the
      same genus if there exists an $\alpha \in k^\times$ such that
      $N \fa = N \fb \cdot N(\alpha)$; note that $N(\alpha) > 0$.
\end{itemize}
As a corollary of genus theory and the index calculations Hecke 
obtains the following characterization: 

\begin{prop}\label{PH}
Let $k$ be a quadratic number field with discriminant $d$. 
An ideal $\fa$ coprime to $d$ is in the principal genus if and 
only if one of the following equivalent conditions is satisfied:
\begin{enumerate}
\item $\fa$ is equivalent in the strict sense to
      the square of some ideal $\fb$; 
\item $(\frac{N\fa,d}{p}) = +1$ for all primes $p \mid d$;
\item $N\fa = N(\alpha)$ for some $\alpha \in k^\times$;
\item $N\fa \equiv N(\alpha) \bmod d$ for some $\alpha \in k^\times$.
\end{enumerate}
\end{prop}

Hecke -- of course -- gives the analytic proof of the existence of
genera:
\begin{quote}
Die Tatsache, da\ss{} die Anzahl der Geschlechter $g$ genau $= 2^{t-1}$
ist, wird nun am bequemsten mit Benutzung transzendenter Methoden 
[\ldots] bewiesen\footnote{The fact that the number $g$ of genera is
$= 2^{t-1}$ can be proved most conveniently by using transcendental
methods [\ldots].}
\end{quote}
and remarks, after having given the `fundamental theorem on genera':
\begin{quote}
Gau\ss{} hat diesen Satz zuerst gefunden und f\"ur ihn einen rein
arithmetischen Beweis gegeben.\footnote{Gau\ss{} has discovered this
theorem and gave a purely arithmetic proof.}
\end{quote}

\medskip
\noindent{\bf Remark.}
Olga Taussky \cite[\S 5]{Tau} characterizes the principal 
genus using matrices.

\section{Euler's Conjectures Revisited}

In this section we will show that Euler's Conjecture 1 and his 
`corrected' Conjecture 2 follow from genus theory. Assume that $n$ 
is a positive squarefree integer and that $p \equiv 1 \bmod 4n$ is 
prime. Then $(p/p_i) = +1$ for all primes $p_i \mid n$, which by
quadratic reciprocity (note that inverting Legendre symbols is no 
problem since $p \equiv 1 \bmod 4$) implies $(d_i/p) = (p_i/p) = +1$, 
where $d_i$ are the prime discriminants\footnote{A prime discriminant 
is a discriminant of a quadratic number field that is a prime power.} 
dividing $d = \disc \Q(\sqrt{n}\,)$. Applying the following 
proposition with $a = -n$ then shows that the Goldbach-Euler
conjecture is true (at least for squarefree $n$):

\begin{prop}\label{PE}
Let $a$ be a squarefree integer $\ne 1$, $k = \Q(\sqrt{a}\,)$
a quadratic number field with discriminant $d$, and $p>0$ a 
prime not dividing $d$. Then the following conditions are
equivalent:
\begin{enumerate}
\item there exist $x, y \in \Q$ with $p = x^2 - ay^2$;
\item we have $(d_i/p) = 1$ for all prime discriminants $d_i$
      dividing $\disc k$;
\item we have $p\OO_k = \fp\fp'$, and $\fp$ is equivalent to the 
      square of an ideal in $\OO_k$ in the strict sense.
\end{enumerate}
\end{prop}

\begin{proof}
Condition (1) says that the norm of a prime ideal $\fp$ above 
$p$ is the norm of an element, which by Proposition \ref{PH}.3
implies that $\fp$ is in the principal genus, which is (3).
Similarly, if $p$ does not divide $d$, then the Legendre
symbols $(d_i/p)$ essentially coincide with Hilbert symbols
$(\frac{p,d}{p_i})$, where $p_i$ is the unique prime dividing 
$d_i$, and this time we see that $\fp$ is in the principal genus 
by Proposition \ref{PH}.2. The claim that (3) $\Lra$ (1) is proved
by simply taking norms.
\end{proof}

Let us now have a look at Lagrange's counterexample to Euler's 
original conjecture in the light of Gau\ss's (or rather Hecke's) 
genus theory. First we observe that $79$ is the smallest natural 
number $a$ such that the class group of $\Q(\sqrt{a}\,)$ is 
strictly larger than the genus class group, and the fact 
that Lagrange didn't give a smaller counterexample suggests 
a connection between Euler's conjecture and genus theory. 

In fact, if we replace the condition $x, y \in \Z$ in
$p = x^2 - Ay^2$ by $x, y \in \Q$, then Euler's
conjecture is true and essentially equivalent to the
principal genus theorem. To the best of my knowledge,
however, this has not been noticed before, neither
by Lagrange (who pointed out that Euler's conjecture
was false), nor by Legendre (who proved a related
result on the solvability of $ax^2 + by^2 + cz^2$,
which contains criteria for the solvability of
$-c = aX^2 + bY^2$ in rational numbers as a special
case), nor by anyone else for that matter. The following
lemma shows the connection of Euler's criterion with
something more familiar: 

\begin{lem}
Let $a$ be a squarefree integer $\ne 1$, $k = \Q(\sqrt{a}\,)$
a quadratic number field with discriminant $d$, and $p>0$ a 
prime not dividing $d$. Then the following conditions are
equivalent:
\begin{enumerate}
\item there exist $n, r \in \Z$ such that $p = 4an + r^2$ or
      $p = 4an+r^2 - a$;
\item we have $(d_i/p) = 1$ for all prime discriminants $d_i$
      dividing $\disc k$.
\end{enumerate}
\end{lem}

The proof is a simple exercise using the quadratic reciprocity law;
Proposition \ref{PE} then shows that Euler's Conjecture 2 is in fact
just a version of the principal genus theorem of quadratic forms or
fields.

\medskip
\noindent{\bf Remark.}
In his preface to Euler's Opera Omnia I--4 Fueter \cite[p. xiii]{Fue} 
remarks that Euler's observation in \cite{EuN} that only half 
of all possible prime residue classes mod $4n$ may yield prime 
factors of $x^2 + ny^2$ is equivalent to Gau\ss{}'s result
that at most half of all possible genera exist. Gau\ss{} had 
already remarked in \cite[\S 151]{Gauss} that there is a gap in
Euler's proof. 

\bigskip 

Edwards \cite[p. 274]{Edw} observes 
 \begin{quote}
   The case $D = 79$ is one that Gau\ss{} frequently uses as an 
example (DA arts 185, 186, 187, 195, 196, 198, 205, 223)
 \end{quote}
and speculates that Gau\ss's interest in this discriminant may have
been sparked by Lagrange's counterexample.

\vskip 1cm

\begin{center} 
{\huge  III. Genus Theory and Higher Reciprocity Laws}
\end{center}

\section{Kummer}

Kummer's motivation for creating a genus theory for Kummer 
extensions over $\Q(\zeta)$, where $\zeta$ is a primitive
$\ell$-th root of unity, was his quest for a proof of 
the reciprocity law for $\ell$-th powers: call an 
$\alpha \in \Z[\zeta]$ primary if $\alpha$ is congruent 
to a nonzero integer modulo $(1-z)^2$ and if 
$\alpha \overline{\alpha}$ is congruent to an integer
modulo $\ell$. For primary and coprime integers
$\alpha, \beta \in \Z[\zeta]$, Kummer had conjectured the
reciprocity law $(\alpha/\beta) = (\beta/\alpha)$, where
$(\,\cdot\,/\,\cdot\,)$ is the $\ell$-th power residue 
symbol. When everything else had failed (in particular 
cyclotomic methods like Gau\ss{} and Jacobi sums), he 
turned to Gau\ss{}'s genus theory.

Let us fix some notation: $\ell$ will denote an odd prime, 
$\zeta$ a primitive $\ell$-th root of unity, $\lambda = 1 - \zeta$, 
and $\fl = (\lambda)$ the prime ideal in $k = \Q(\zeta)$ above 
$\ell$. Let $M$ denote the set of all $\alpha \in k^\times$ 
coprime to $\fl$.

The generalization of Gau\ss's theory of genera of quadratic forms
to Kummer extensions $k(\sqrt[\ell]{\alpha}\,)/k$ looks quite
natural to mathematicians who are familiar with the interpretation
of the theory of binary quadratic forms as an arithmetic of ideals
in quadratic number fields; although Kummer never worked out such
a theory, he was aware that his ideal numbers could be generalized 
to this situation (with the same problems as over Kummer extensions 
of cyclotomic fields: primes dividing the conductor of the ring had 
to be excluded).

Armed with this insight, it was clear that the quadratic extensions
$\Q(\sqrt{d}\,)$ had to be replaced by Kummer extensions, and that
representability of primes by binary quadratic forms corresponded
to being norms of prime ideals in ideal classes.\footnote{We 
replace Kummer's language by Dedekind's. We also use the notation
familiar from Galois theory; instead of denoting the conjugates of
$\alpha = f(\zeta)$ by $f(\zeta^r)$ as Kummer did, we let 
$\sigma: \zeta \too \zeta^r$ act on $\alpha$.} 
But what are the right analogs of Gau\ss's characters? As Kummer 
eventually realized, these characters could be constructed using 
`differential logarithms'.

Assume that $\alpha \in \Z[\zeta]$ satisfies 
$\alpha \equiv 1  \bmod \lambda$, and write it as 
$\alpha = f(\zeta)$ for some polynomial $f \in \Z[X]$; 
replace $X$ by the function $e^v$, evaluate the $r$-th 
derivative of $\log f(e^v)$ with respect to $v$ at $v=0$, 
and call the result $\dl^r(\alpha)$ (Kummer wrote 
$\frac{d_0^r \log f(e^v)}{dv^r}$ instead; the subscript $0$ 
at $d_0^r$ indicates that the derivative should be evaluated 
at $v = 0$). For $1 \le r \le \ell-2$, the resulting integer 
modulo $\ell$ does not depend on the choice of $f$; with a 
little bit more care it can be shown that a similar procedure 
gives a well defined result even for $r = \ell-1$.

It was noticed by Takagi and Hasse that Kummer's differential logarithms
can (and should) be described algebraically; here is a short summary
of the most basic properties of the $\dl^r$ (see \cite[Chapter 14]{L2}): 

\begin{prop}\label{PKd}
Kummer's differential logarithms\footnote{I prefer differential 
logarithm to logarithmic derivative since \ref{PKd}.i) resembles 
the functional equation of the logarithm and not the Leibniz rule.} 
$\dl^r$ $(1 \le r \le p-1)$ have the following properties:
\begin{enumerate}
\item[i)] $\dl^r(\alpha\beta) = \dl^r(\alpha) + \dl^r(\beta)$ 
          for all $\alpha, \beta \in M$
\item[ii)] $\dl^r(\alpha^\sigma) = \chi(\sigma)^r\dl^r(\alpha)$,
           where $\chi: \Gal(k/\Q) \lra \F_\ell$ is the cyclotomic
           character defined by $\zeta^\sigma = \zeta^{\chi(\sigma)}$;
\item[iii)] $\dl^r(\alpha) = \dl^r(\beta)$ for all $\alpha, \beta \in M$
          with $\alpha \equiv \beta \bmod \lambda^{r+1}$;
\item[iv)] $\dl^r(1+\lambda^r) = (-1)^r r!$.
\end{enumerate}
\end{prop}

Kummer introduced his differential logarithms in \cite[p. 493]{Kum}
in connection with Gauss sums and observes in a special case that 
$\dl^r(\alpha^\sigma) = (-1)^r \dl^r(\alpha)$, where $\sigma$ 
denotes complex conjugation. This special case of Proposition
\ref{PKd}.ii) immediately implies that $\dl^r(\alpha) = 0$ for 
all real $\alpha \in k$ whenever $r$ is odd; in particular,
$\dl^{2r+1}(\eps) = 0$ for all real units in $\Z[\zeta]$.

Put $K = \Q(\zeta_\ell)$, fix an integer $\mu \in \Z[\zeta_\ell]$
and consider the Kummer extension $L = K(\srt{\ell}{\mu})$.
Kummer's ``integers in $w$'' were elements of
$\OO[w]$, where $w = \srt{\ell}{\mu}$ and $\OO = \Z[\zeta_\ell]$;
observe that $\OO[w] \ne \OO_L$ in general even when $\mu$ is
squarefree.  Next he introduces integers
$z_j = (1-\zeta)(1-\mu)/(1-w\zeta^j) \in \OO[w]$ as well as 
the ring $\OO_z = \OO[z_0, z_1, \ldots, z_{\ell-1}]$ and 
observes that $\OO_z \subseteq \ell\OO[w] \subseteq \OO_z$. 
On [p. 676-677] he discusses the decomposition law, excluding 
primes dividing $\lambda\mu$.

Assume that $\fp$ is a prime ideal in $\Z[\zeta]$ and
let $h$ denote the class number. Then $\fp^h = (\pi)$,
and we can try to define $\dl^r(\fp)$ by the equation
$h \dl^r(\fp) = \dl^r(\pi)$. Unfortunately, the values 
of $\dl^r(\pi)$ depend on the choice of $\pi$ in general, 
although not always: since $\dl^{2r+1}(\eps_j) = 0$ for all 
real units $\eps_j$ and all $0 \le r \le \rho = \frac12(\ell-1)$, 
and since moreover $\dl^{2r+1}(\zeta) = 0$ for all $1 \le r \le \rho$, 
we can define $\dl^{2r+1}(\fp) = \frac1h \dl^{2r+1}(\pi)$.

This allows Kummer to define $\rho$ characters $\chi_3$, $\chi_5$, 
\ldots, $\chi_{\ell-2}$ on the group of ideals in $\OO_z$ prime to 
$\ell \cdot \disc(L/K)$ by putting
$$\chi_{2r+1}(\fP) = \zeta^S, \quad S = \dl^{2r+1}(N_{L/K}\fP).$$
To these characters he adds
$$\chi_{\ell-1}(\fP) = \zeta^S, \quad S = \frac{1-N\fP}{\ell}.$$

Now let $\fp_1, \ldots, \fp_t$ denote the primes different from 
$(\lambda)$ that are ramified in $L/K$. For each such prime Kummer 
defines a character 
$$\psi_j(\fP) = \bigg(\frac{N_{L/K}\fP}{\fp_j}\bigg).$$
Here $\fp = N_{L/K}\fP$ is an ideal in $\Z[\zeta]$,
$\fp^h = (\pi)$ is principal, and if we insist in taking
$\pi$ primary, then the symbol $(\pi/\fp_j)$ only depends
on $\fP$. We therefore put $(N_{L/K}\fP/\fp_j) = (\pi/\fp_j)^{h^*}$,
where $h^*$ is an integer such that $h^*h \equiv 1 \bmod \ell$. 

All in all there are now $\rho + t$ characters, and these can be 
shown to depend only on the ideal class of $\fP$ (\cite[p. 748]{Kum}).
The ideal classes with trivial characters form a subgroup
$C^z_\gen$ in $\Cl^z(L)$, the class group of the order
$\OO_z$, and $C^z_\gen$ is called the principal genus.
The quotient group $\Cl^z_\gen(L/K) = \Cl^z(L)/C^z_\gen$
is called the genus class group, and the main problem is
to determine its order. This problem is solved by invoking
ambiguous ideal classes (see \cite[p. 751]{Kum}): 
\begin{quote}
Die Anzahl aller wirklich vorhandenen Gattungen ist nicht
gr\"o\ss{}er, als die Anzahl aller wesentlich verschiedenen,
nicht \"aquivalenten ambigen Klassen.\footnote{The number of 
existing genera is not greater than the number of all 
essentially different nonequivalent ambiguous classes.} 
\end{quote}
 
On the next 40 pages, Kummer \cite[p. 752--796]{Kum} shows that 
there are exactly $\ell^{\rho+t-1}$ ambiguous ideal classes; this 
is quite a surprising result, because the ambiguous class number 
formulas the we are familiar with all contain a unit index as
a factor. The amazing thing is that it is Kummer's `weird'
choice of the order he is working in eliminates this index!
By working in an order with a nontrivial conductor Kummer
is actually able to simplify genus theory considerably.

As the number of pages he spends on this topic shows, he had
to work hard nonetheless. For counting the number of ambiguous
ideal classes he comes up with his integers in $u$, coinciding
with the integers in $w$ in the case where $\mu$ is a prime ideal
power. Actually what he is doing is writing the $h$-th power of 
$\mu$ as a product of $t$ principal prime powers, adjoining 
their $\ell$-th roots $u_j$ to $\OO_z$, and then showing that 
ambiguous prime ideals in $\OO_z$ become principal in the 
extension $\OO_u = \OO_w[u_1, \ldots, u_t]$ in such a way that
their generators can be made into elements of $\OO_w$ through 
multiplication by powers of the $u_j$ (\cite[p. 768]{Kum}).

To illustrate this result, take $k = \Q(\sqrt{15}\,)$ with the
ambiguous ideals $\ft = (2,1+\sqrt{15}\,)$, $\fr = (3,\sqrt{15}\,)$
and $\fv = (5,\sqrt{15}\,)$: these become principal in the ring
$\Z[\sqrt{3},\sqrt{5}\,]$ as $\ft = (\sqrt{3} + \sqrt{5}\,)$,
$\fr = (\sqrt{3}\,)$ and $\fv = (\sqrt{5}\,)$, and multiplying 
these generators with $\sqrt{3}$ produces integers in $\Z[\sqrt{15}\,]$.
Had Kummer demanded that the $u_j$ be primary, the $u_j$ would 
generate a subfield of the genus class field of $K/k$, and his 
result would say that ambiguous ideal classes become principal 
there: this is a near miss that even Hilbert did not follow out
(although it may have inspired his Satz 94 on the capitulation
of ideals in unramified cyclic extensions).

Next Kummer \cite[p. 796]{Kum} obtains the first inequality of genus
theory, namely the fact that there are at most $\ell^{\rho+t-1}$ genera 
(Kummer writes $\lambda$ instead of $\ell$):
\begin{quote}
Die Anzahl der wirklich vorhandenen Gattungen ist nicht gr\"o\ss{}er
als der $\ell$-te Theil aller blo\ss{} angebbaren Gattungen oder
Gesamtcharaktere.\footnote{The number of existing genera is not 
greater than the $\ell$-th part of all possible genera or
total characters.}
\end{quote}
He notes, however, that this is not good enough to prove the 
reciprocity law: imitating Gau\ss{}'s second proof only gives 
a distinction between $\ell$-th power residues and nonresidues, 
but is not powerful enough to distinguish between e.g. 
$(\alpha/\fp)_\ell = \zeta$ and $(\alpha/\fp)_\ell = \zeta^2$.

Kummer closes this gap by proving the second inequality
in some special cases. To this end, he has to study norm
residues modulo powers of $(1-\zeta)$ in the Kummer extensions
$\Q(\zeta,\sqrt[\ell]{\mu}\,)/\Q(\zeta)$. His first result 
is that if a number $\alpha \in \Z[\zeta]$ is a norm from 
$\OO_w$, then (\cite[p. 805]{Kum})
\begin{equation}\label{KuNR}
\dl^1(\alpha)\dl^{\ell-1}(\mu) + \dl^2(\alpha)\dl^{\ell-2}(\mu) + 
  \dl^{\ell-1}(\alpha)\dl^1(\mu)  \equiv 0 \bmod \ell.
\end{equation}
This is an amazing result: the left hand side of (\ref{KuNR}) 
is an element of $\F_p$, and this element vanishes if $\alpha$
is a norm from $\OO_w$. Hilbert later realized that the left
hand side is just the additively written norm residue symbol
at the prime $\fp$ above $p$.

On p. 808 Kummer shows that condition (\ref{KuNR}) is equivalent to 
$$ \Big(\frac{\eps}{\mu}\Big) = \Big(\frac{\eta}{\alpha}\Big), $$ 
where $\eps$ and $\eta$ are units such that $\eps\alpha$ and
$\eta\mu$ are primary.

If $\fp$ is a prime ideal in $\OO_k$ such that $(\eps/\fp) = 1$
for all units $\eps \in \OO_k^\times$, then $\fp$ is called
a prime ideal of the second kind, and of the first kind otherwise.
 
The first special case is obtained on p. 811: if $t = 1$, and if
the ramified prime ideal is of the first kind, then there are 
exactly $\ell^\rho$ genera. On p. 817, he obtains a similar
result for certain Kummer extensions with exactly two ramified 
primes. This turns out to be sufficient for proving the reciprocity 
laws, but before he does so, he applies these reciprocity laws to
derive the general principal genus theorem on p. 825:
\begin{quote}
Die Anzahl der wirklich vorhandenen Gattungen der idealen Zahlen 
in $z$ ist genau gleich dem $\ell$-ten Theile aller Gesamtcharaktere.
\footnote{The number of existing genera in the theory of ideal
numbers in $z$ is equal to the $\ell$-th part of all total characters.}
\end{quote}

\section{Hilbert}

Hilbert's Zahlbericht\footnote{It was almost immediately translated 
into French (1909); meanwhile, there also exist translations into
Romanian (1997) and English (1998).} consists of five parts: 
the foundations of ideal theory, Galois theory, quadratic number
fields, cyclotomic fields, and Kummer extensions, and the first four
parts are still considered to be standard topics in any introduction
to algebraic number theory. The fifth part, clearly the most difficult
section of the Zahlbericht, did not make it into any textbook and was 
soon superseded by the work of Furtw\"angler and Takagi. Yet it is 
this chapter that I regard to be the Zahlbericht's main claim to fame: 
it reflects Hilbert's struggle with digesting Kummer's work, of 
finding a good definition of the norm residue symbol, and of 
incorporating Kummer's isolated results on genus theory of Kummer 
extensions into a theory which is on a par with the genus theory
of binary quadratic forms in Gau\ss's Section V of the Disquisitiones.

\medskip

The quadratic norm residue symbol $(\frac{n\,,\,m}{p})$ 
is defined to be $+1$ if $m$ is a square or if $n$ is congruent
modulo any power of $p$ to the norm of a suitable integer from 
$\Q(\sqrt{m}\,)$, and  $(\frac{n\,,\,m}{p}) = -1$ otherwise.
This Hilbert symbol can be expressed using Legendre symbols;
in \cite[Satz 13]{Hil2} Hilbert derives the formula
$$  \Big(\frac{\nu,\mu}{\fp}\Big) = 
                \Big(\frac{(-1)^{ab}\rho\sigma}{\fp}\Big) $$
for the Hilbert symbol for primes $\fp \nmid 2$ in number 
fields $k$, where $\fp^a \parallel \mu$, $\fp^b \parallel \nu$,
and $\nu^a\mu^{-b} = \rho\sigma^{-1}$ for integers 
$\rho, \sigma \in \OO_k$ coprime to $\fp$.

For defining the $\ell$-th power norm residue symbol for odd primes 
$\ell$, Hilbert proceeds in the opposite direction. Let $\vc_\fp$ 
denote the $\fp$-adic valuation, i.e.' let $\vc_\fp(\alpha)$ be
the maximal power of $\fp$ dividing $\alpha$. For
$\mu, \nu \in k^\times$ put $a = \vc_\fp(\mu)$
and $b= \vc_\fp(\nu)$. Then $\nu^a\mu^{-b} = \rho\sigma^{-1}$
for integers $\rho,\sigma \in \OO_k$ such that
$\vc_\fp(\rho) = \vc_\fp(\sigma) = 0$ . Now define
$$  \Big(\frac{\nu,\mu}{\fp}\Big)_\ell^{\phantom{p}} =
                \Big(\frac{\rho}{\fp}\Big)_\ell^{\phantom{p}}
                \Big(\frac{\sigma}{\fp}\Big)_\ell^{-1} $$
for all prime ideals $\fp$ not dividing $\ell$.

The definition of the norm residue symbol for primes ideals
$\fp \mid \ell$ is much more involved; in his Zahlbericht,
Hilbert only considers the case $k = \Q(\zeta_\ell)$ and uses
Kummer's differential logarithms in the case $\ell \ge 3$:
for $\mu \equiv \nu \equiv 1 \bmod \ell$, he puts (compare
Kummer's result (\ref{KuNR}))
$$  \Big(\frac{\nu,\mu}{\fl}\Big)_\ell^{\phantom{p}} =
  \zeta^S \ \  \text{with} \
 S \ = \ \dl^1(\nu)\dl^{\ell-1}(\mu) - \dl^2(\nu)\dl^{\ell-2}(\mu)  \pm
   \ldots - \dl^{\ell-1}(\nu)\dl^1(\mu),$$
and then extends it to $\mu, \nu$ coprime to $\ell$ by
$$  \Big(\frac{\nu,\mu}{\fl}\Big)_\ell^{\phantom{p}} \ = \
   \Big(\frac{\nu^{\ell-1},\mu^{\ell-1}}{\fl}\Big)_\ell^{\phantom{p}}. $$

Hilbert's genus theory goes like this:\footnote{We have taken the 
liberty of rewriting it slightly using the concept of quotient groups.} 
let $k = \Q(\zeta_\ell)$, and assume that the class number $h$ of $k$ 
is not divisible by $\ell$. Consider the Kummer extension
$K = k(\sqrt[l]{\mu}\,)$.  Let $\fp_1$, \ldots, $\fp_t$
denote the primes that are ramified in $K/k$ (including
infinite ramified primes if $\ell = 2$). For each ideal
$\fa$ in $\OO_k$, write $N_{K/k}\fa^h = \alpha\OO_k$; the map
$$   \alpha \too X(\alpha) =
   \Big\{ \Big(\frac{\alpha,\mu}{\fp_1} \Big), \ldots,
          \Big(\frac{\alpha,\mu}{\fp_t} \Big) \Big\} $$
induces a homomorphism $\psi: \Cl(K) \lra {\F_l}^t/X(E_k)$
by mapping an ideal class $[\fa]$ to $X(\alpha)^{h^*}X(E_k)$,
where $h^*$ is an integer with $h^*h \equiv 1 \bmod \ell$.
Its kernel $C_\gen = \ker \psi$ is called the principal
genus, and the quotient group $\Cl_\gen(K) = \Cl(K)/C_\gen$
the genus class group of $K$. 

In \cite[Satz 150]{Hil} Hilbert generalizes Gau\ss's observation 
on p. \pageref{GO} by proving that the index of norm residues modulo 
$\fp^e$ in the group of all numbers coprime to $\fp$ is $1$ if 
$\fp$ is unramified, and equal to $\ell$ if $\fp \ne \fl$ is 
ramified or if $\fp = \fl$ and $e > \ell$.

Next Hilbert shows that his symbol defined in terms of power residue
symbols actually is a norm residue symbol in \cite[Satz 151]{Hil}.

Following Gau\ss{}, Hilbert then derives the inequality $g \le a$ 
between the number of genera and ambiguous ideal classes 
(\cite[Hilfssatz 34]{Hil}), then proves the reciprocity law 
$\prod_v (\frac{\,a\, , \, b\,}{v}) = 1$ for the $\ell$-th 
power Hilbert symbol and regular primes $\ell$ (\cite[\S 160]{Hil}), 
and finally obtains the second inequality $g \ge a$ 
(\cite[Satz 164]{Hil}). This result is then used for proving the
principal genus theorem in \cite[Satz 166]{Hil}:
\begin{quote}
jede Klasse des Hauptgeschlechtes in einem regul\"aren Kummerschen
K\"orper $K$ ist gleich dem Produkt aus der $1-S$-ten sym\-bo\-li\-schen
Potenz einer Klasse und einer solchen Klasse, welche Ideale des
Kreisk\"orpers $k(\zeta)$ enth\"alt.\footnote{every class of the 
principal genus in a regular Kummer field $K$ is the product of  
the $1-S$-th symbolic power of an ideal class and of a class 
containing ideals of the cyclotomic field $k(\zeta)$.}
\end{quote}
(Observe that this implies the familiar equality 
$C_\gen = \Cl(K)^{1-\sigma}$ if we work with $\ell$-class groups.)
Satz 167 finally shows that numbers in $k$ that are norm residues at 
every prime $\fp$ actually are norms from $K$, and Hilbert concludes
this section with the remark
\begin{quote}
Damit ist es dann gelungen, alle diejenigen Eigenschaften auf den
regul\"aren Kummerschen K\"orper zu \"ubertragen, welche f\"ur den
qua\-dra\-ti\-schen K\"orper bereits von {\sc Gauss} aufgestellt und
bewiesen worden sind.\footnote{Thus we have succeeded in transferring
all those properties to the regular Kummer fields that already have 
been stated and proved for quadratic number fields by Gau\ss{}.}
\end{quote}

For connections between genus theory and reciprocity laws see also 
Skolem \cite{Skol}.

\vskip 1cm

\begin{center} 
{\huge  IV. Genus Theory in Class Field Theory}
\end{center}

\section{Furtw\"angler}

In Furtw\"angler's construction of Hilbert class fields,
the following theorem (see \cite{Fur07}) played a major role: 

\begin{thm}\label{TFw}
Let $L/K$ be a cyclic unramified extension, $\sigma$ a generator
of the Galois group $\Gal(L/K)$, and let $N: \Cl(L) \lra \Cl(K)$
be the norm map on the ideal class groups. Then 
$\ker N = \Cl(L)^{1-\sigma}$. 
\end{thm}

This almost looks like the vanishing of $H^{-1}(G,\Cl(L))$, but 
actually we have $H^{-1}(G,\Cl(L)) \ne 0$ in general; this is due 
to the difference between the relative norm 
$$N_{L/K}: \Cl(L) \lra \Cl(K)$$ 
and the algebraic norm 
$$\nu_{L/K} = 1 + \sigma + \sigma^2 + \ldots + \sigma^{(L:K)-1}: 
        \Cl(L) \lra \Cl(L);$$ 
the connection between these two norms is the relation
$\nu = j \circ N$, where $j: \Cl(K) \lra \Cl(L)$ is the transfer
of ideal classes. This means that Furtw\"angler's principal
theorem can't be translated easily into the cohomological 
language because ideal classes may capitulate. Furtw\"angler 
\cite{Furt} used his principal genus theorem to study the 
capitulation of ideals in Hilbert $2$-class fields of 
number fields with $2$-class group isomorphic to $(2,2)$. 

Furtw\"angler also proved that, for cyclic extensions $L/K$
of prime degree, an element $\alpha \in K^\times$ is a norm 
from $L$ if and only if it is a norm residue modulo the 
conductor $\frf$ of $L/K$ (Hasse's contribution was the
interpretation of this result as a Local-Global principle). 
We will later see that this result can be expressed 
cohomologically as $H^{-1}(\Gal(L/K),C_L) = 1$, where 
$C_L$ is the id\`ele class group; for this reason, Kubota 
\cite{Kub} calls $H^1(C) = 0$ the principal genus theorem 
and credits Furtw\"angler for the `fully id\`ele-theoretic' 
result in the case of Kummer extensions of prime 
degree.\footnote{In \cite{Kub}, Kubota shows that the 
second inequality of class field theory is essentially
a corollary of two of Furtw\"angler's results: the product
formula for the Hilbert symbol (i.e., the reciprocity law),
and the principal genus theorem mentioned above.}
Nakayama \cite{Nak}, on the other hand, claims that $H^1(C) = 0$ 
is `merely the id\`ele-class analogue of Noether's principal genus 
theorem', while Chevalley \cite{Chev} calls it the generalization 
of Hasse's principal genus theorem to normal extensions.

\section{Takagi and Hasse}

In this section we assume some familiarity with the
classical version of class field theory. 
Let $L/K$ be an extension of number fields and $\fm$ a
modulus in $K$. Let $P^1\{\fm\}$ denote the set of
principal ideals $(\alpha)$ in $K$ with $\alpha \equiv
1 \bmod \fm$, and let $D_K\{\fm\}$ denote the group of
ideals in $K$ coprime to $\fm$, and let $D_K\{\fm\}$
denote the corresponding object for $K$. Then we call
$H_{L/K}\{\fm\} = N_{L/K} D_L\{\fm\} \cdot P^1\{\fm\}$
the ideal group defined mod $\fm$ associated to $L/K$.

In the special case where $\fm$ is an integral ideal,
such groups had been studied by Weber; in their theory
of the Hilbert class field, Hilbert and Furtw\"angler 
defined infinite primes, and Takagi combined these two
notions to create his class field theory.

Takagi called $L$ a class field of $K$ for the ideal
group $H_{L/K}\{\fm\}$ if $(D_K\{\fm\} : H_{L/K}\{\fm\})
 = (L:K)$. In order to show that abelian extensions are
class fields, this equality has to be proved, and the 
proof is done in two steps: 
\begin{enumerate}
\item the First Inequality 
      $$ (D_K\{\fm\}:H_{L/K}\{\fm\}) \le (L:K),$$
      holds for any finite extension $L/K$ and any modulus $\fm$
      and can be proved rather easily using analytic techniques.
\item the Second Inequality says that 
      $$ (D_K\{\frf\}:H_{L/K}\{\frf\}) \ge (L:K),$$
      where $L/K$ is a cyclic extension of prime degree $l$, and
      where $\frf$ is the conductor of $L/K$, that is, the ideal
      such that the relative discriminant of $L/K$ is $\frf^{l-1}$.
\end{enumerate}
In his famous Marburg lectures \cite{HasM} on class field theory,
Hasse puts the proof of the second inequality into historical
perspective by mentioning the role of Gau\ss's work:
\begin{quote}
Wir gehen jetzt auf den Beweis des Umkehrsatzes aus. Die dazu
erforderlichen \"Uberlegungen des laufenden Paragraphen bilden 
die Verallgemeinerung der ber\"uhmten Gau\ss{}schen Untersuchungen
\"uber die Theorie der Geschlechter quadratischer Formen aus seinen
{\em disquisitiones arithmeticae}.\footnote{We now are going for the
proof of the inverse theorem. The considerations of this section,
which  will be needed to achieve this, are generalizations of 
Gau\ss's famous investigations in the genus theory of quadratic
forms in his {\em disquisitiones arithmeticae}.}
\end{quote}

In the following, we will provide the background needed
for portraying the role of genus theory in the proof of 
the second inequality of class field theory.

Consider a  cyclic extension $L/K$ of prime degree $\ell$;
then $\disc(L/K) = \frf^{\ell-1}$ for some ideal $\frf$ in 
$\OO_K$ called the conductor of $L/K$. Takagi's definition 
of genera in $L/K$ is based on the observation that there 
is a connection between the class group $\Cl(L)$ and some
ray class group $\Cl_K^\nu$ defined modulo $\frf$: given a 
class $c = [\fA] \in \Cl(L)$, we can form the ray class 
$[N_{L/K} \fA]$ in the group $\Cl_K^\nu$ of ideals modulo norm 
residues, that is, in the group $D_K\{\frf\}$ of ideals coprime 
to $\frf$ modulo the group $P_K^\nu\{\frf\}$ of principal ideals 
generated by norm residues modulo the conductor $\frf$. Note that 
if $\fA = \lambda \fB$ for some $\lambda \in L^\times$, then 
the ray classes generated by $N_{L/k} \fA$ and $N_{L/k} \fB$ 
coincide since $N_{L/K} \lambda \in P_K^\nu\{\frf\}$.

Consider the norm map $N_{L/K}: \Cl(L) \lra \Cl_K^\nu$.
Clearly $\im N = N \Cl(L) = H_{L/K}\{\frf\}/P_K^\nu\{\frf\}$,
so the image of the norm involves the ideal group
associated with $L/K$. The kernel of the norm map
is called the principal genus $C_\gen$: it is the group 
of all ideal classes $c = [\fA] \in \Cl(L)$ such that 
$N_{L/K}\fA = (\alpha)$ for norm residues $\alpha \in K^\times$ 
(thus $\alpha$ is coprime to $\frf$ and a norm residue at 
every prime ideal). Thus we find the following exact sequence
\begin{equation}\label{ESG}  \begin{CD} 
  1 @>>> C_\gen @>>> \Cl(L) @>{N}>> H_{L/K}\{\frf\}/P_K^\nu\{\frf\} @>>> 1. 
   \end{CD} \end{equation}
Thus computing the number of genera $g = (\Cl(L):C_\gen)$ will 
help us in getting information about the order of the ideal class 
group associated to $L/K$. We will show that $g = a$, where $a$ 
denotes the number of ambiguous ideal classes in $L$. In fact,
$C_\gen$ clearly contains the group $\Cl(L)^{1-\sigma}$, where 
$\sigma$ is a generator of $\Gal(L/K)$. This shows that
$$a = (\Cl(L):\Cl(L)^{1-\sigma}) \ge (\Cl(L):C_\gen) = g,$$
that is, the first inequality of genus theory.
The left hand side can be evaluated explicitly: the ambiguous
class number formula says that
$$ a = h_K \cdot \frac{\prod e(\fp)}{(L:K)(E:E_\nu)}, $$
where $h_K = \# \Cl(K)$ is the class number of $K$, 
$e(\fp)$ is the ramification index of a prime ideal $\fp$ in $L/K$,
the product is over all (ramified) primes in $K$ including the
infinite primes, $E$ is the unit group of $K$, and $E_\nu$ its 
subgroup of units that are norm residues modulo $\frf$.
  
For proving the second inequality of genus theory, namely $g \ge a$, 
we use the exact sequence (\ref{ESG}) and get 
$$ (\Cl(L):C_\gen)  = (N \Cl(L):1) = (H_{L/K}\{\frf\}:P_K^\nu\{\frf\}) 
   = \frac{(D_K\{\frf\}: P_K^\nu\{\frf\})}{(D_K\{\frf\}:H_{L/K}\{\frf\})}.$$
Now the index in the denominator satisfies
$(D_K\{\frf\}:H_{L/K}\{\frf\}) \le \ell$ by the first inequality.
The index in the numerator is the product of 
$h_K = (D_K\{\frf\}:P_K\{\frf\})$ (ideals away from $\frf$ modulo
the subgroup of principal ideals), that is, the class number of $K$, 
and the index $(P_K\{\frf\}:P_K^\nu\{\frf\})$; this last index can be 
computed explicitly, and it turns out that
$$ (P_K\{\frf\}:P_K^\nu\{\frf\}) = (E_\nu:E \cap NL^\times) 
              \cdot \frac{\prod e(\fp)}{(E:E_\nu)}. $$
Thus we find
$$ (D_K\{\frf\}: P_K^\nu\{\frf\}) =
             (E_\nu:E \cap NL^\times) \cdot a\ell \ge a\ell;$$
this means that in the sequence of inequalities
$$ a  \ge (\Cl(L):C_\gen) = g = 
      \frac{(D_K\{\frf\}: P_K^\nu\{\frf\})}
           {(D_K\{\frf\}:H_{L/K}\{\frf\})} \ge a$$
we must have equality throughout; in particular we find
\begin{itemize}
\item $(D_K\{\frf\}:H_{L/K}\{\frf\}) = \ell:$ cyclic extensions 
      are class fields;
\item The principal genus theorem: $C_\gen = \Cl(L)^{1-\sigma}$; note
      that if $L/K$ is unramified, then $C_\gen = \Cl(L)[N]$ coincides
      with the kernel of the norm map $\Cl(L) \lra \Cl(K)$, and the
      principal genus theorem becomes Theorem \ref{TFw}.  
\item The norm theorem for units: $(E_\nu:E \cap NL^\times) = 1$, 
      that is, any unit that is a norm residue modulo the conductor 
      is the norm of some element of $L^\times$.
\end{itemize}
Thus the proof of the second inequality consists in a 
calculation of the number $g$ of genera: the inequality 
$g \le a$ comes from the ambiguous class number formula,
the inequality $g \ge a$ from the first inequality of
class field theory and some cohomological results.

Takagi then derives the norm theorem (in cyclic extensions, norm 
residues modulo the conductor are actual norms) from the
principal genus theorem.

\subsection*{Hasse} 
In his Klassenk\"orperbericht \cite{Has}, Hasse reproduces 
Takagi's proof of the second inequality with only minor
modifications. In his Marburg lectures \cite{HasM}, on the 
other hand, Hasse proves the second inequality 
\begin{equation}\label{EHa}
   (D_K\{\frf\}:H_{L/K}\{\frf\}) \ge (L:K)
\end{equation}
in a different and direct way; the main advantages of his proof are 
\begin{itemize}
\item it is valid for cyclic extensions of arbitrary (finite) degree;
\item the full norm theorem is a consequence of equality in (\ref{EHa});
\item it does not use the first inequality.
\end{itemize}
This last fact later allowed Chevalley to give an arithmetic proof 
of class field theory by proving the second inequality first and 
then deriving the first inequality without analytic means.

At some point in the computation of (\ref{EHa}), the index
(norm residues modulo conductor : norms) is written as the 
product of (units that are norm residues : norms of units) and 
(ideal classes of the principal genus : $(1-\sigma)$-th powers 
 of ideal classes). Thus Hasse's norm theorem (which follows by
comparing (\ref{EHa}) with the first inequality) contains the 
principal genus theorem, i.e., the statement that the principal 
genus consists of the $(1-\sigma)$-th powers of ideal classes.

\subsection*{The General Principal Genus Theorem}
Recall that ideal classes in $L$ were mapped by the norm to
ray classes modulo $\frf$ in $K$. Are there similar results
connecting ray classes in $L$ with ray classes in $K$? The 
answer is yes: in \cite[p. 304--310]{Has}, Hasse proved the 
`most general' principal ideal theorem. In order to state it 
we need the following 

\begin{prop}
Let $L/K$ be a cyclic extension of prime degree with generating
automorphism $\sigma$, and let $\fm$ be a modulus in $K$. Then
there exists a modulus $\fM$ in $L$ such that
\begin{enumerate}
\item $\fM \mid \fm \OO_K$;
\item $\fM^\sigma = \fM$;
\item for $\beta \in L^\times$ coprime to $\fM$ we have 
      $N_{L/K}(\beta) \equiv 1 \bmod \fm$  if and only if 
      $\beta \equiv \alpha^{1-\sigma} \bmod \fM$. 
\end{enumerate}
\end{prop}

With these preparations, Hasse defines the principal genus
$\bH_1$ mod $\fM$ in $L$ as the group of ray classes modulo $\fM$
whose relative norms land in the ray modulo $\fm$ in $K$. 

\begin{thm}\label{HPG}
Let $L/K$, $\fm$ and $\fM$ be defined as above. Then the
principal genus $\bH_1$ coincides with the group of 
$1-\sigma$-th powers of ray classes mod $\fM$ in $L$. 
\end{thm}

This was generalized even more by Herbrand \cite{Herb}.

\section{Chebotarev and Scholz}

The generalization of genus theory from cyclic extensions
to general normal extensions was mainly the work of
Chebotarev \cite{Cheb} and Scholz \cite{Sch}.

Let $L/K$ be a normal extension. The maximal unramified 
extension of $L$ of the form $LF$, where $F/K$ is abelian,
is called the genus class field $L_\gen$ of $L$ with respect to $K$; 
the maximal unramified extension that is central over $K$ is called 
the central class field and is denoted by $L_\cen$. 

According to Scholz, these definitions are due to Chebotarev \cite{Cheb};
as a matter of fact, his paper is not easy reading. The characterization
of the genus and central class fields in terms of class groups is
due to Scholz:

\begin{thm}\label{TSch}
Let $L/K$ be a normal extension of number fields, let $H_0$ 
denote the elements of $K^\times$ that are norm residues, 
and put $N_0 = N_{L/K} L^\times$. Next, let $H$ and $N$ denote 
the group of ideals in $L$ whose norms land in the groups of
principal ideals generated by elements of $H_0$ and $N_0$, 
respectively.\footnote{Observe that $H$ is the principal genus 
in the sense of Takagi.} Then the class field associated to the 
ideal group $H$ is the genus class field $L_\gen$, and the class
field associated to $N$ is the central class field $L_\cen$. 
In particular, Scholz's number knot $H_0/N_0$ is isomorphic to
the Galois group of $L_\cen/L_\gen$. 
\end{thm}

Scholz used this result to prove that Hasse's norm residue theorem 
(everywhere local norms are global norms) is valid in all extensions
whose Galois groups have trivial Schur multiplier. Jehne \cite{J}
presented Scholz's work in a modern language and extended his
results. 

As an unramified abelian extension of $L$, $L_\gen$ corresponds to
some quotient $\Cl(K)/C_\gen$ of the class group of $K$, and the
group $C_\gen$ is called the principal genus, which for cyclic
extensions $L/K$ can be shown to satisfy $C_\gen = \Cl(L)^{1-\sigma}$, 
where $\sigma$ is a generator of $\Gal(L/K)$.

The following theorem connects the modern definition of the
principal genus with the classical one by Takagi: Fr\"ohlich 
\cite[pp. 18--19]{FroC} calls it the {\em classical principal 
genus theorem}:

\begin{thm}
Let $L/K$ be a cyclic extension, and $\sigma$ a generator of $\Gal(L/K)$. 
Then $[\fa] \in C_\gen$ if and only if $N_{L/K} \fa = (\alpha)$, where 
$\alpha \in K^\times$ is a norm residue at all ramified primes in $L/K$. 
\end{thm}

This form of genus theory was used by various number theorists;
among the many contributions, let us mention Hasse \cite{HasG}
and Leopoldt \cite{Leo}, Gold \cite{Gold}, Stark \cite{Sta} 
(whose generalization of genus theory lacks an analogue of 
Gau\ss's principal genus theorem), Gurak \cite{Gur}, and Razar 
\cite{Raz}. 

\vskip 1cm

\begin{center} 
{\huge  V. Genus Theory and Galois Cohomology}
\end{center}

\section{Noether}

Emmy Noether thought very highly of her version of the principal 
genus theorem for number fields that she developed in early 1932.
It was published in 1933 in a paper \cite{HGS} which became more 
famous for the `Noether equations'\footnote{Lorenz \cite{Lo} has 
observed that these are due to Speiser, and that Noether actually
credits him in \cite{HGS}.} in connection with Hilbert's Theorem 90
than for the main content, the principal genus theorem.
Expositions of Noether's principal genus theorem can be 
found in Deuring \cite[VII, \S 7]{Deu} and Fr\"ohlich \cite{Froe}.

\subsection*{Noether's version}
Noether starts with a short introduction to crossed products:
let $K$ be a field and $L/K$ a separable extension of degree $n$
and with Galois group $G$. The crossed product of $L$ and $G$
is an algebra $A$ together with injections $L \hookrightarrow A$
and $G \hookrightarrow A$ such that all automorphisms of $L$
become inner automorphisms of $A$.

Noether next describes this algebra $A$ using factor systems. 
As a $L$-vector space, $A$ is generated by the basis elements
$u_{\sigma_1}$, \ldots, $u_{\sigma_n}$ corresponding to the $n$ group
elements $s_i$; thus we have
\begin{equation}\label{EN1}
A = u_{\sigma_1} L \oplus \cdots \oplus u_{\sigma_n}L.
\end{equation}
The condition that the automorphisms $\sigma$ on $L$ should become inner 
can be satisfied by demanding
\begin{equation}\label{EN2}
  z^\sigma = u_\sigma^{-1} z u_{\sigma}
\end{equation}
for every $z \in L$. This defines a factor system $(a_{\sigma,\tau})$ 
in $L^\times$ by
\begin{equation}\label{EN3}
  u_\sigma u_\tau = s_{\sigma\tau} a_{\sigma,\tau},
\end{equation}
and associativity of multiplication gives the relation
\begin{equation}\label{EN4}
 a_{\sigma\tau,\rho} a_{\sigma,\tau}^\rho 
               a_{\sigma,\tau\rho} a_{\tau,\rho}.
\end{equation}
Now the multiplication
$$ \sum u_\sigma b_\sigma \cdot \sum u_\tau c_\tau 
             = \sum u_\sigma u_\tau b_\sigma^\tau c_\tau  $$
makes $A$ into a simple normal algebra over $K$ which will be 
denoted by $A = (a_{\sigma,\tau}, L, G)$. Different factor systems 
$a_{\sigma,\tau}$ and $\ba_{\sigma,\tau}$ generate isomorphic 
algebras if there are $c_\sigma \in L^\times$ such that
\begin{equation}\label{EN5}
  \ba_{\sigma,\tau} = a_{\sigma,\tau} c_\sigma^\tau c_\tau / c_{\sigma\tau}.  
\end{equation}
The cosets $u_\sigma L^\times$ define a group extension $G^\times$ of $G$:
$$ \begin{CD} 1 @>>> L^\times @>>> G^\times @>>> G @>>> 1. \end{CD} $$

The principal genus theorem Noether is about to prove has an analogue
for normal extensions, namely Hilbert's Theorem 90: she gives three
different formulations of this result.

\begin{prop}\label{PNM}
Let $L/K$ be a finite Galois extension with Galois group $G$.
Then the following assertions are equivalent formulations of the
`Minimal' Principal Genus Theorem: 
\begin{enumerate}
\item Every group automorphism of $G^\times$ whose restriction 
      to $L^\times$ is the identity is inner, and is generated
      by an element of $L^\times$.
\item If $c_\sigma^\tau c_\tau / c_{\sigma\tau} = 1$ 
      for all $\sigma, \tau \in G$, then there exists a 
      $b \in L^\times$ such that $c_\sigma = b^{1-\sigma}$ 
      for all $\sigma \in G$.
\item The group $G$ has a unique crossed representation class of
      first degree associated to the trivial factor system.
\end{enumerate}
\end{prop}

Using the language of cohomology groups, the second 
version of the minimal principal genus theorem claims 
that $H^1(G,L^\times) = 0$, i.e., it is Hilbert's Satz 
90. A representation $u_\sigma \too C_\sigma$ is called 
a crossed representation for the factor system 
$a_{\sigma,\tau}$ if 
$C_\sigma^\tau C_\tau = C_{\sigma\tau} a_{\sigma,\tau}$. 
Two crossed representations $u_\sigma \too C_\sigma$ and 
$u_\sigma \too D_\sigma$ for $a_{\sigma,\tau}$ belong to 
the same class if $C_\sigma = B^{-\sigma} D_\sigma B$.

Noether then defines an {\em ideal factor system} of 
a Galois extension $L/K$ with Galois group 
$\Gal(L/K) = \{\sigma, \tau, \ldots\}$ as a system of $n^2$ 
ideals $\fa_{\sigma,\tau}$ of $L$ satisfying the relations 
$$ \fa_{\sigma,\tau} \fa_{\tau, \rho} =
     \fa_{\sigma\tau,\rho} \fa_{\sigma,\tau}^\rho.$$
>From $n$ ideals $\fc_\sigma$ we can construct a factor system
$$ \fa_{\sigma,\tau} = \frac{\fc_\sigma^\tau \fc_\tau}{\fc_{\sigma\tau}}$$
called the transformation system.

Ideal factor systems form a group $C$, and the transformation systems
form a subgroup $B$ of $C$.  In analogy to the group of norm residues 
modulo the conductor Noether defines the principal class of ideal
factor systems as consisting of systems $\fa_{\sigma,\tau}$ with 
the following property: there exists a factor system $a_{\sigma,\tau}$ 
in $L^\times$ such that
\begin{enumerate}
\item $\fa_{\sigma,\tau} = (a_{\sigma,\tau})$;
\item $a_{\sigma,\tau}$ determines an algebra $\fA = (L,a)$ that
      splits at every ramified place $\fp$ of $L/K$.
\end{enumerate}

Then Noether's principal genus theorem states

\begin{thm}
If the transformation system 
$\fc_\sigma^\tau \fc_\tau\fc_{\sigma\tau}^{-1}$ is in the
principal class, then there is an ideal class $[\fb]$ such that
$[\fc_\sigma] = [\fb]^{1-\sigma}$ for all $\sigma \in \Gal(L/K)$. 
\end{thm}

Noether also gives two formulations analogous to the first and
third version of Proposition \ref{PNM}. 

\subsection*{Id\`eles}
For stating Noether's principal genus theorem in a modern language
we need to introduce the the id\`ele group $J$ of a number field $L$, 
the id\`ele class group $C$, the unit id\'eles $U$, the group of 
fractional ideals $D$ and its subgroup $P$ of principal ideals, as 
well as other well known invariants of $L$. These groups are all 
part of the fundamental exact and commutative square
$$ \begin{CD} 
   @. 1 @. 1 @. 1 @. \\
   @. @VVV @VVV @VVV @. \\
 1 @>>> E @>>> L^\times @>>>  P @>>> 1 \\
   @. @VVV @VVV @VVV @. \\
 1 @>>> U @>>> J @>>> D @>>> 1 \\
   @. @VVV @VVV @VVV @. \\
1 @>>> {\mathcal E} @>>> C @>>> \Cl @>>> 1 \\
   @. @VVV @VVV @VVV @. \\
   @. 1 @. 1 @. 1 @.   \end{CD} $$
Our aim is to reformulate Noether's principal genus theorem 
in the language of the cohomology of id\`eles; fix a normal 
extension $L/K$ of number fields with Galois group $G = \Gal(L/K)$. 
As all our cohomology groups will be formed with $G$, we put 
$H^q(M) := H^q(G,M)$ from now on. 

There is a `Hilbert 90' for id\`eles: we have $H^1(J_L) = 1$;
this is essentially a direct consequence of Hilbert's Theorem 
90 for the localizations $L_\fp$. The claim $H^1(C_L) = 1$ is much 
deeper; in fact, it can be viewed as a generalization of Hasse's 
norm theorem for cyclic extensions, that is, as a Local-Global 
Principle for normal extensions.

In order to see this, take the long cohomology sequence of the
exact sequence
$$ \begin{CD} 1 @>>> L^\times @>>> J_L @>>> C_L @>>> 1; \end{CD} $$
observing that $H^{-1}(J_L) \simeq H^1(J_L) = 1$ for cyclic
groups $G$ we find
$$ \begin{CD} 
   1 @>>> H^{-1}(C_L) @>>> H^0(L^\times) @>i>> H^0(J_L).
   \end{CD} $$
Now $H^0(L^\times) = K^\times/N_{L/K} L^\times$, $H^0(J_L) = J_K/N J_L$, 
hence $\ker i = \partial_{L/K} := K^\times \cap N J_L/NL^\times$ is
the obstruction to Hasse's Local-Global Principle. Thus
$H^1(C) = 1$ is equivalent to Hasse's norm theorem. Observe that
$\partial_{L/K} = H_0/N_0$ coincides with Scholz's number knot
from Theorem \ref{TSch}.

\subsection*{Roquette's version}
Let us now translate Noether's results into the language of 
cohomology of id\`eles. We will start with a version due 
to P.~Roquette \cite{Ro} and then explain how this is connected 
with A.~Fr\"ohlich's `translation'.
 
Let $H^2_S(L^\times)$ denote the subgroup of $H^2(L^\times)$
whose elements split at all primes in $S$; in other words,
$H^2_S(L^\times)$ is the kernel of the natural map
$H^2(L^\times) \lra H^2\big(\prod_{w \in S} L_w^\times\big)$
induced by $L^\times \lra \prod_{w \in S} L_w^\times$
(where $w \in S$ is short for $w \mid v$ for $v \in S$).

Next the exact sequence $1 \lra E \lra L^\times \lra P \lra 1$ 
yields a map $H^2_S(L^\times) \lra H^2(P)$; let $H_S$ denote the 
image of $H_S^2(L^\times)$ in $H^2(P)$.

Finally, the sequence $1 \lra P \lra D \lra \Cl \lra 1$ 
gives us an exact sequence 
\begin{equation}\label{N1}  \begin{CD} 
   0 = H^1(D) @>>> H^1(\Cl) @>{\delta}>> H^2(P) @>{\phi}>> H^2(D),
\end{CD} \end{equation}
which allows us to identify $H^1(\Cl)$ with a subgroup of $H^2(P)$.
In particular, we can talk about $H_S \cap H^1(\Cl)$.

Using this language, Noether's principal genus theorem can be 
stated in the following way:

\begin{thm}\label{TPGN}
Let $L/K$ be a normal extension of number fields with Galois 
group $G$, and let $S$ be the set of ramified primes. 
Then $H_S \cap H^1(\Cl) = 0$.
\end{thm}

In order to see the conection with Noether's original formulation,
observe that a transformation system $\fc_\sigma$ of ideals is a 
cocycle of the ideal class group and therefore defines an element 
$c_\sigma \in H^1(\Cl)$ if 
$[\fc_\sigma]^\tau [\fc_\tau] = [\fc_{\sigma\tau}]$ (this is 
the first condition of the system $\fc_\sigma$ being in the principal 
class); the claim that $c_\sigma = [\fc_\sigma] = [\fb]^{1-\sigma}$ 
says that the cocycle is a coboundary, in other words, that 
$H^1(\Cl) = 1$. This is, however, only true if the second condition 
is also satisfied; this condition demands that 
$\fc_\sigma^\tau \fc_\tau\fc_{\sigma\tau}^{-1} = (a_{\sigma,\tau})$
for a factor system $a_{\sigma,\tau} \in H^2(L^\times)$ whose
associated algebra splits at every place $\fp$ that is ramified 
in $L/K$. In our language, the element $c_\sigma \in H^1(\Cl)$ 
defines a factor set of principal ideals in $H^2(P)$ under the 
connection homomorphism; if this factor system actually comes 
from an element in $a_{\sigma,\tau} \in H^2(L^\times)$ whose 
associated algebra splits at the ramified primes (i.e., if 
$a_{\sigma,\tau} \in H_S^2(L^\times)$), then Noether's principal 
genus theorem claims that the element $c_\sigma$ is trivial. 

\subsection*{Proof of Noether's Theorem}
Let us now derive Noether's principal genus theorem from 
class field theory. An ingredient of Noether's  proof of 
the principal genus theorem is what she calls the 

\begin{prop}[Principal Genus Theorem for Ideals]
If $L/K$ is a normal extension of number fields with Galois group $G$,
then $H^1(D_L) = 0$, where $D_L$ is the group of fractional ideals in $L$.
\end{prop}

Our first goal is to prove 

\begin{lem}\label{Linj}
We have $H^1(C) = 0$ if and only if $H^2_S(L^\times) \lra H^2(D)$
is injective. 
\end{lem}

\begin{proof}
Since $H^1(J) = 0$, we find that the Local-Global Principle 
$H^1(C) = 0$ holds if and only if $H^2(L^\times) \lra H^2(J)$ is
injective. 

>From the definition of $H_S^2(L^\times)$ and the analogous
definition of $H_S^2(J)$ we get the exact and commutative
diagram
$$ \begin{CD}
   0 @>>> H_S^2(L^\times) @>>> H^2(L^\times) 
             @>>> H^2\big(\prod_{w \in S} L_w^\times\big) \\
   @.  @.  @VVV @VVV \\
   0 @>>> H_S^2(J) @>>> H^2(J) 
             @>>> H^2\big(\prod_{w \in S} L_w^\times\big). \end{CD}$$
A simple diagram chase shows that the map  $H^2(L^\times) \lra H^2(J)$
induces a map $H_S^2(L^\times) \lra H_S^2(J)$, and that 
$\ker (H^2(L^\times) \lra H^2(J)) = \ker (H^2_S(L^\times) \lra H^2_S(J))$. 

The natural map $J \lra D$ from id\`eles to their ideals induces 
a homomorphism $H^2(J) \lra H^2(D)$, so by restriction we get a
map $H^2_S(J) \lra H^2(D)$. 

We have already seen that $H^1(C) = 0$ if and only if 
$H^2(L^\times) \lra H^2(J)$ is injective, and that 
$\ker (H^2(L^\times) \lra H^2(J)) =  
 \ker (H^2_S(L^\times) \lra H^2_S(J))$.

Thus it remains to show that $H^2_S(J) \lra H^2(D)$ is injective.
The exact sequence $1 \lra U \lra J \lra D \lra 1$ gives, via
the long cohomology sequence, rise to a commutative diagram
$$ \begin{CD}
  0 = H^1(D) @>>> H^2(U) @>>> H^2(J) @>>> H^2(D) \\
   @.  @VVV @VVV @VVV \\
   0   @>>> H^2\big(\prod_{w \in S} L_w^\times\big) 
       @>>> H^2\big(\prod_{w \in S} L_w^\times\big) @>>> 0. \end{CD} $$
Applying the snake lemma gives an exact sequence 
$$ \begin{CD} 0 @>>> H_S^2(U) @>>> H_S^2(J) @>>> H^2(D) \end{CD} $$
It remains to show that $H_S^2(U) = 0$ if $S$ contains ramified primes;
observe that $U \simeq \prod_{w \in S} U_w \times \prod_{w \notin S} U_w$,
hence $H^2(U) \simeq H^2(\prod_{w \in S} U_w) \times
       H^2(\prod_{w \notin S} U_w)$.   

It is known (see e.g. \cite[p. 131, Prop. 1.1]{CF}) that
$H^q(\Gal(L/K), U_L) = 0$ for $q \ge 1$ and unramified extensions
$L/K$ of local fields; globally we have 
$H^q(G, \prod_{w \mid v} U_w) \simeq H^q(G_{w_0},U_{w_0})$
by \cite[p. 176/177, Prop. 7.2]{CF}, where $w_0$ is a fixed 
prime in $L$ above $v$, and in particular we have 
$H^2(U) = H^2(\prod_{w \mid v \in S} U_w)$, hence
$H_S^2(U) = 0$ as claimed.

Thus $\ker (H^2_S(L^\times) \lra H^2_S(J)) = 
 \ker (H^2_S(L^\times) \lra H^2(D))$ as claimed
\end{proof}

Now we can prove Noether's Principal Genus Theorem \ref{TPGN}:
 
\begin{proof}
Take an element $\delta(c) \in H^1(\Cl)$; if $\delta(c) \in H_S$,
then there is an element $a \in H_S^2(L^\times)$ such that
$\delta(c) = \alpha(a)$, where $\alpha: H_S^2(L^\times) \lra H^2(P)$.
Since (\ref{N1}) is exact, we have $0 = \phi \circ \delta(c)
= \phi \circ \alpha (a)$. Since $\phi \circ \alpha$ is injective
by Lemma \ref{Linj}, we conclude that $a = 0$, then $0 = \alpha(a)
= \delta(c)$ implies $c = 0$ since $\delta$ is injective by 
(\ref{N1}).
\end{proof}

\medskip
\noindent{\bf Remark.}
If $L/K$ is cyclic and if $\sigma$ generates $\Gal(L/K)$, then
$H_S \cap H^1(\Cl) \simeq H_S \cap H^{-1}(\Cl)$, where the intersection
is taken in $H^0(P)$. We find that $H_S$ is the image of 
$H^0_S(L^\times)$ in $H^0(P)$, i.e., that $H_S$ is the subgroup 
of all principal ideals generated by norm residues at the primes in 
$S$ modulo principal ideals generated by norms.    

On the other hand, $H^{-1}(\Cl) = \Cl(L)[\nu]/\Cl(L)^{1-\sigma}$,
where $\nu = 1 + \sigma + \ldots + \sigma^{(L:K)-1}$ is the
algebraic norm $\nu: \Cl(L) \lra \Cl(L)$. The image of 
$H^{-1}(\Cl)$ in $H^0(P) = P^G/P^\nu$ under the connecting 
homomorphism is computed as follows: take $c = [\fa]$ in a
class representing an element of $H^{-1}(\Cl)$; then 
$\fa^\nu = (A)$ for some $A \in L^\times$, and clearly
$(A)^\sigma = (A)$, so $(A) \in P^G$. The class of $c$ is
then mapped to the class of $(A)$. 

This means that $H_S \cap H^{-1}(\Cl)$ consists of cosets of 
ideal classes $c \in \Cl(L)$ such that $c = [\fa]$ with
$\fa^\nu = (\alpha)$, where $\alpha \in K^\times$ 
is a local norm at the primes in $S$. In other words: If
$L/K$ is cyclic and $S$ contains the ramified primes, then
$H_S \cap H^{-1}(\Cl) = C_\gen/\Cl(L)^{1-\sigma}$, and 
Noether's theorem implies that $C_\gen = \Cl(L)^{1-\sigma}$.
If $L/K$ is unramified, then we can take $S = \varnothing$ 
and get back Furtw\"angler's Theorem \ref{TFw}.

\subsection*{Fr\"ohlich's version}
Fr\"ohlich's article \cite{Froe} contains a cohomological 
interpretation of Noether's principal genus theorem that differs 
slightly from Roquette's version above. Here's Fr\"ohlich's account: 
Let $J'$ denote the ideles that have entries $1$ at all places
outside of $S$; then $J' \simeq \prod_{w \in S} L^\times_w$.
The projection $J \lra J'$ and the inclusion $L^\times \lra J$
give rise to a maps $\pi_1: H^2(J) \lra H^2(J')$
and with $\iota: H^2(L^\times) \lra H^2(J)$; thus
$\ker \pi_1 \circ \iota = H^2_S(L^\times)$. 

Next he defines maps $\psi: H^2(L^\times) \lra H^2(P)$
and $\phi: H^2(P) \lra H^2(D)$ (he uses $I$ instead of $D$);
then $\psi (\ker \pi_1 \circ \iota) = H_S$. 

The equivalent of the injectivity of $H_S^2(L^\times) \lra H^2(D)$
is the statement that 
$\ker \pi_1 \circ \iota \cap \ker \phi \circ \psi = 0$: in fact,
$\ker \pi_1 \circ \iota = H^2_S(L^\times)$ and 
$\ker \phi \circ \psi = \ker (H^2(L) \lra H^2(D))$. 

Fr\"ohlich's version of Noether's theorem reads:

\begin{thm}
The composition of maps 
$$H^1(\Cl) \lra H^2(P) \lra H^2(P)/\psi (\ker \pi_1 \circ \iota)$$
is injective.
\end{thm}

In our language: the induced map $H^1(\Cl) \lra H^2(P)/H_S$ is
injective. This is, of course, equivalent to the statement
that $H_S \cap H^1(\Cl) = 0$ in $H^2(P)$.



\subsection*{Terada, Tannaka, Deuring}

Noether's formulation of the principal genus theorem apparently
wasn't very influential, but there have been a few articles
picking up her ideas.

Terada \cite{Ter1} stated Hasse's principal genus theorem in
the following form: given a cyclic extension $L/K$, a generator
$\sigma$ of $\Gal(L/K)$, and modules $\fm$ and $\fM$ as in 
Theorem \ref{HPG}, the following assertions are equivalent:
\begin{eqnarray}
\label{Te1} N_{L/K} \fA \equiv 1 \bmod \fm 
        & \text{for some ideal}\ \fA \subseteq \OO_L; \\
\label{Te2} \fA = \fB^{1-\sigma} \alpha 
        & \text{with}\ \alpha \equiv 1 \bmod \fM. 
\end{eqnarray} 

He introduces the crossed homomorphism $\fA_\tau: G \lra D_L$
by putting $\fA_1 = (1)$, $\fA_\sigma = \fA$, \ldots, and
$\fA_{\sigma^a} = \fA_\sigma^{\sigma^{a-1}} \fA_{\sigma^{a-1}}$
for $0 < a < (L:K)$. Then he observes that (\ref{Te1}) is equivalent
to the condition 
\begin{equation}\label{Te3}
   \fA_\rho^\tau \fA_\tau \fA_{\rho\tau}^{-1} \equiv 1 \bmod \fm
\end{equation}
for all $\rho, \tau \in G$, and that (\ref{Te2}) is equivalent to the
existence of an ideal $\fB$ such that
\begin{equation}\label{Te4}
 \fA_\rho^\tau = \fB^{1-\sigma} \alpha_\rho, 
          \quad  \alpha_\rho \equiv 1 \bmod \fM.
\end{equation}
These numbers $\alpha_\rho$ satisfy the conditions
\begin{equation}\label{Te5} 
\alpha_\rho^\tau \alpha_\tau \alpha_{\rho\tau}^{-1} \equiv 1 \bmod \fm
\end{equation}
for all  $\rho, \tau \in G$. 

Terada's generalized version of the principal genus theorem then is

\begin{thm}
Let $L/K$ be an abelian extension; then a system $\fA_\rho$ satisfies the
condition (\ref{Te3}) if and only if there exist an ideal $\fB$ in $L$ 
and elements $\alpha_\rho \in L^\times$ such that (\ref{Te4}) and 
(\ref{Te5}) are satisfied.
\end{thm}

For related versions, see Kuniyoshi \&  Takahashi \cite{KT}
and Terada \cite{Ter2}. The last result we will mention is 
a conjecture of Deuring proved by Tannaka \cite{Tan}:

\begin{thm}
Let $L$ be the Hilbert class field of a number field $K$. For every ideal
$\fa$ in $\OO_K$ there is an element $\Theta(\fa) \in L^\times$ with the
following properties:
\begin{enumerate}
\item $\fa\OO_L = \Theta(\fa) \OO_L$;
\item Let $\sigma(\fa) = (\frac{L/K}{\fa})$ be the Artin automorphism, and define
the factor set $\eps(\fa,\fb)$ by
$$ \eps(\fa,\fb) = 
   \frac{\Theta(\fa) \Theta(\fb)^{\sigma(\fa)}}{\Theta(\fa\fb)}. $$
Then $\eps(\fa,\fb) \in E_K$, the unit group of $\OO_K$.  
\end{enumerate} 
\end{thm}

\section*{From Gau\ss\ to Noether}

We have seen how the principal genus theorem
\begin{itemize}
\item evolved from a striking result in the theory of binary 
      quadratic forms, 
\item was translated into the language of quadratic number fields 
      and generalized in Kummer's and Hilbert's quest for general
      reciprocity laws, 
\item became the main technique for proving the second inequality 
      in Takagi's class field theory, 
\item and finally ended up as a collection of lemmas about the 
      vanishing of (parts of) the first cohomology of certain 
      groups associated to number fields. 
\end{itemize}
In a somewhat parallel development, the formulation of 
reciprocity laws (which motivated the development of genus 
theory) changed considerably: there is
\begin{itemize}
\item the classical version of Euler, Legendre, and Gau\ss;
\item the product formula for the Hilbert symbol;
\item the isomorphism of Artin's reciprocity law;
\item the bijection in the (still largely hypothetical)
      formulation within the Langlands correspondence.
\end{itemize}
The big difference, however, is this: whereas the reciprocity 
law has maintained its central place in abelian as well as nonabelian 
class field theory, genus theory went from the center of interest 
into almost complete oblivion. It is even stranger that its climax 
(the work of Hasse, Herbrand, and Noether in the 1930s) and its 
sudden death happened almost simultaneously. 

\section*{From Noether to Gau\ss}

In a letter\footnote{See HINT \cite{HINT}.} to Hasse from 
June 2/3, 1932, Noether writes:
\begin{quote}
Im \"ubrigen habe ich anl\"a\ss{}lich der Ausarbeitung
meines Z\"uricher Vortrags einmal Gau\ss{} gelesen. Es
wurde behauptet, da\ss{} ein halbwegs gebildeter
Mathematiker den Gau\ss{}schen Hauptgeschlechtssatz
kennt, aber nur Ausnahmemenschen den der
Klassenk\"orpertheorie. Ob das stimmt, wei\ss{} ich nicht
-- meine Kenntnisse gingen in um\-ge\-kehr\-ter Reihenfolge --
aber jedenfalls habe ich in bezug auf Auffassung allerhand
von Gau\ss{} gelernt; vor allem da\ss{} es gut ist den
Nachweis, da\ss{} die durch Faktorensysteme bestimmte
Klasseneinteilung eine Strahlkl.-Einteilung ist, an den
Schlu\ss{} zu stellen; der \"Ubergang meiner Fassung zu
der Gau\ss{}schen geht n\"amlich unabh\"angig davon direkt,
erst die Spezialisierung auf die Klassenk\"orpertheorie
braucht den F\"uhrer. Was ich mache, ist die
Verallgemeinerung der Definition der Geschlechter durch
Charaktere.\footnote{By the way, during the preparation 
for my Zurich lecture I read Gau\ss{}. It has been
claimed that a reasonably educated mathematician knows
Gau\ss{}'s principal genus theorem, but only exceptional
people the principal genus theorem of class field theory.
I don't know if that's true -- in my case it was the other
way round -- but in any case I learned a lot about
perception from Gau\ss{}; above all that it is a good idea
to place the proof that the classes determined by factor 
systems are ray classes at the end; the transition from
my version to Gau\ss{}'s can be done independently and directly,
only for the specialization to class field theory the conductor
is needed. What I am doing is the generalization of the
definition of genera using characters.}
\end{quote}
 
The paragraph on conductors was omitted from the final version:
in another letter to Hasse on October 29, 1932, Noether sends
him a final draft of \cite{HGS} and writes
\begin{quote}
Den F\"uhrer-Paragraphen habe ich fortgelassen; er wurde so
kompliziert da\ss{} es keinen Sinn hatte, da nichts damit
gemacht wurde, \ldots\footnote{I have omitted the paragraph
on conductors; it became so complicated that it made no sense
[to keep it], since I did not use it, \ldots}
\end{quote}
On November 25, 1932, Noether thanks Hasse for his comments on
her article and writes
\begin{quote}
Den F\"uhrer habe ich weggelassen; das wird einmal zusammen mit
anderem kommen, wenn es mehr durchgearbeitet ist.\footnote{I omitted
the conductor; this will be taken care of together with other
things once I have worked it over.}
\end{quote}

With this bow of Emmy Noether to Gau\ss{} we conclude our
survey of the development of the principal genus theorem.

\section*{Acknowledgements}
I would like to thank Peter Roquette for explaining the content
of Noether's paper to me, as well as Norbert Schappacher for
numerous helpful comments.

\end{document}